\documentclass[12pt]{article}
\usepackage[centertags,intlimits]{amsmath}                     
\usepackage{color}\usepackage{amsfonts}                     
\usepackage{amsthm}
\usepackage{graphics}
\usepackage{graphicx}
\allowdisplaybreaks[2]
\numberwithin{equation}{section}
\newtheorem{theorem}{Theorem}[section]
\newtheorem{lemma}{Lemma}[section]
\newtheorem{corollary}{Corollary}[section]

\newtheorem{remark}{Remark}[section]

\providecommand{\abs}[1]{\lvert #1\rvert}

\newcommand{\nc}{\newcommand}
\nc{\vb}{\mathbf{v}}
\nc{\bx}{\mathbf{x}}
\nc{\by}{\mathbf{y}}
\nc{\bz}{\mathbf{z}}
\nc{\bu}{\mathbf{u}}
\nc{\bv}{\mathbf{v}}
\nc{\ba}{\mathbf{a}}
\nc{\bs}{\mathbf{s}}
\nc{\bq}{\mathbf{q}}
\nc{\bd}{\mathbf{d}}
\nc{\bb}{\mathbf{b}}
\nc{\bc}{\mathbf{c}}
\nc{\bi}{\mathbf{i}}
\nc{\bfr}{\mathbf{r}}
\nc{\bA}{\mathbf{A}}
\nc{\R}{\mathbb R}
\nc{\N}{\mathbb N}
\nc{\D}{\mathbb D}
\nc{\Z}{\mathbb Z}
\nc{\F}{\mathbf F}
\nc{\bbS}{\mathbb S}
\nc{\B}{\cal B}
\nc{\br}{\bigr}
\nc{\bl}{\bigl}
\nc{\Bl}{\Bigl}
\nc{\Br}{\Bigr}
\nc{\ind}{\mathbf{1}}
\nc{\bP}{\mathbf{P}}

\title{Large deviations, optimal trajectories and conservation of momenta
in generalised Jackson networks}
\author{Anatolii A. Puhalskii \footnote{Email: puhalski@iitp.ru}\\
 Institute for Problems in Information
Transmission}
\begin{document}

\maketitle
\sloppy
\begin{abstract} 
We study properties of optimal trajectories for
 the quasipotential that governs large
deviations of queue lengths in generalised Jackson networks (GJN). 
The trajectories are shown to
satisfy Hamilton equations. 
The
Hamiltonian is identified.
Optimal trajectories and  associated costs
are found for networks with two and three stations.
A conservation law for adjoint variables
  is discovered to the effect that certain components are
maintained along an optimal trajectory. 
A process  is proposed
for finding optimal trajectories and computing their costs
 in  general GJN.
\end{abstract}

\section{Introduction }
\label{sec:introduction}

The variational principle for large deviations identifies
the   probability of the rare event that a
 stochastic process 
 visits   a neighbourhood of a certain point $x$ in the
 state space over large time  with $e^{-nV(x)}$\,,  $n$ being a large
 parameter
  and $V(x)$ representing the quasipotential.
A rigourous formulation requires  a  Large Deviation
Principle
(LDP) for the
stochastic process and
asserts asymptotic logarithmic equivalence for the probability in question,
 see Freidlin and Wentzell \cite{wf2},
 Shwartz and Weiss \cite{ShwWei95}.
 Justifications of
the principle are confined for the most part to the Markovian
framework. In Puhalskii \cite{Puh19a,Puh_25b} a version of the variational
principle was established for the non Markov process of 
queue lengths in an ergodic
generalised Jackson network (GJN).
In this paper we address the issue of evaluating the quasipotential.
As the qiasipotential 
 is the minimum of the values of the trajectorial deviation
function, a.k.a. rate function, over trajectories that pass through $x$,
we look for
optimal trajectories. Those
trajectories are also representative of most likely scenarios that realise
large queue lengths, see, e.g. Shwartz and Weiss \cite{ShwWei95},
Puhalskii \cite{Puh19}.

Similarly to Puhalskii \cite{Puh19}, the
trajectories are  piecewise linear, being linear in the regions of
constant dynamic, i.e. in the regions with the same number of zero coordinates,
which  we call facets, and
have slope changes when transiting between facets.
The trajectories solve Hamilton equations,  the
Hamiltonian being associated with the trajectorial deviation function in
Puhalskii \cite{Puh07}. We identify the Hamiltonian  for a
general GJN.
For two and three station networks,
we find  trajectories that get optimally from the origin to 
points of the state space. For the two station network a  graphical
illustration is provided.
 Adjoint variables,  which we call
momenta in analogy with classical mechanics,
  play a primary role in our study,  rather than
coordinates. We discover a conservation law for the momenta to the effect
that momentum components that correspond to positive coordinates are maintained along an optimal trajectory.
 We
conjecture that this property holds generally and propose a process
for finding  momenta in GJN.

The paper is organised as follows. 
Section \ref{sec:quasp-hamilt-equat}
defines the model,
 introduces  notation and  assumed hypotheses, provides
necessary background and states the main result.
 Sections \ref{sec:two-stations} and 
\ref{sec:three-stations} analyse networks of two and three stations,
respectively. Section \ref{sec:some-generalities} concerns
generalisations and conjectures.
Proofs and other technical details are collected in the appendices.
\section{The quasipotential and Hamiltonian}
\label{sec:quasp-hamilt-equat}

We consider a  queueing network with
 a homogeneous customer population which comprises $K$   
 single server stations.   
Customers arrive exogenously  and are served  in the
order of arrival, one customer at a time. Upon being served, they either 
join a queue at another station or leave the network.
Let $A_k(t)$  denote the   
cumulative number of exogenous arrivals  at station $k$ 
by time $t$\,,
let $S_k(t)$ denote the cumulative number of customers   
that are served   
at station $k$ for the first $t$ units of busy time of that
station, and let 
$\Phi_{kl}(m)$
denote the cumulative number of customers among the   
first $m$ customers  departing station $k$ that go directly to station   
$l$. Let
$A_k=(A_k(t),\,t\in\R_+)$,
$S_k=(S_k(t),\,t\in\R_+)$, and $\Phi_k=(\Phi_{k}(m),\,m\in\Z_+)$,
where $\Phi_{k}(m)=(\Phi_{kl}(m),\,l\in\mathcal{K})$ and
 $\mathcal{K}=\{1,2,\ldots, K\}$\,. 
It is assumed that the $A_k$ and $S_k$ are nonzero ordinary renewal processes
and
$  \Phi_{kl}(m)=\sum_{i=1}^m  1_{\{\zeta_{k}^{(i)}=l\}}$,
where $\{\zeta_k^{(1)},\zeta_k^{(2)},\ldots\}$ is  a sequence of
i.i.d. r.v. assuming values in  $\mathcal{K}\cup\{0\}$\,,
$ 1_\Gamma$ standing for the indicator function of set $\Gamma$\,.
 Let $Q=(Q(t),\,t\in\R_+)$ represent 
  the queue length process,
where   $Q(t)=(Q_k(t),\,k\in\mathcal{K})$ and 
$Q_{k}(t)$ represents the number of customers  at station $k$ at time   
$t$\,.  The random entities $A_k$\,, $S_l$\,, $\Phi_m$ and $Q(0)$ are
assumed
to be defined on a common probability
space $(\Omega,\mathcal{F},\mathbf{P})$ and
be mutually independent,  where 
$k,l,m\in\mathcal{K}$\,.
All  the stochastic processes are assumed  to have
piecewise constant right--continuous with
left--hand limits trajectories. Accordingly, they are considered as random
elements of the associated Skorohod spaces.

We  denote $p_{kl}=\mathbf P(\zeta_k^{(1)}=l)$ and let $P=(p_{kl})_{k,l=1}^K$\,.
The matrix $P$ is assumed to be of spectral radius less than unity so
that every arriving customer eventually leaves.
For $k\in\mathcal{K}$   and $t\in\R_+$, 
 the following equations hold:
\begin{equation}
  \label{eq:2}
Q_k(t)=Q_k(0)+A_k(t)+\sum_{l=1}^K\Phi_{lk}\bl(D_{l}(t)\br)-D_k(t),
\end{equation}
where 
\begin{equation}
  \label{eq:1}
D_k(t)=S_k\bl(B_k(t)\br)
\end{equation}
represents  the number of  departures  from station $k$ by time $t$ and
\begin{equation}
  \label{eq:5}
  B_k(t)=\int_0^t 1_{\{Q_k(u)>0\}}\,du
\end{equation}
represents the cumulative busy time of station $k$ by time $t$\,.
For  given realisations of
 $Q_k(0)$, $A_k$,  $S_k$, and $\Phi_k$, 
there exist unique  $Q_k=(Q_k(t),\,t\in\R_+)$, 
$D_k=(D_k(t),\,t\in\R_+)$ and $B_k=(B_k(t),\,t\in\R_+)$ 
that satisfy 
 \eqref{eq:2}, \eqref{eq:1} and \eqref{eq:5}, see, e.g.,
Chen and Mandelbaum \cite{MR1106809}.

Let nonnegative random variables $\xi_k$ and
$\eta_k$ represent generic times between exogenous arrivals  and
service times at station $k$, respectively. It is assumed that
$0<\mathbf E\xi_k<\infty$ and  $0<\mathbf E\eta_k<\infty$\,.
Let  $\lambda_k=1/\mathbf E\xi_k$ and $\mu_k=1/\mathbf E\eta_k$\,.
We  assume that the network  is subcritical, i.e.
\[
\mu>(I-P^T)^{-1}\lambda
\]
where 
$\mu=(\mu_1,\ldots,\mu_K)^T$\,,
$\lambda=(\lambda_1,\ldots,\lambda_K)^T$ and
 $I$ represents the $K\times K$--identity matrix.
(Vector inequalities are understood  entrywise.)
In addition, it is assumed that the interarrival times are
unbounded and spreadout.
 More specifically, it is assumed that
 \begin{enumerate}
\item $\mathbf P(\xi_k>u)>0$\,, for all
$k\in\mathcal{K}$ and  $u>0$\,,
\item for each $k\in\mathcal{K}$\,,
there exist nonnegative function $\phi_k(u)$ on $\R_+$ with
  $\int_0^\infty \phi_k(u)\,du>0$ and $m_k\in\N$ 
such that $\mathbf P(\sum_{i=1}^{m_k}\xi_{k,i}\in[x,y])\ge
\int_x^y\phi_k(u)\,du$\,, provided  $0\le x\le y$\,, where 
$\xi_{k,1},\ldots,\xi_{k,m_k}$ are i.i.d. and are distributed as $\xi_k$\,.
 \end{enumerate}
Let $U(t)=(U_1(t)\,, \ldots,U_K(t))$ represent the $K$--vector of excess times of
the exogenous arrival processes at $t$\,. To put it another way, 
$U_k(t)$ is the
length of time that elapses
 between $t$  and the next exogenous arrival after $t$ at station
$k$\,.
 Let $W(t)=(W_1(t)\,, \ldots,W_K(t))$ represent the
$K$--vector of residual service times  at the
stations  at $t$\,, i.e. $W_k(t)$ is the length of time until the first
service completion after $t$\,,
with the understanding that $W_k(t)=0$ if station $k$ is idle at $t$\,.
Under the
introduced hypotheses, 
   $X(t)=(Q(t),U(t),W(t))$ is a Harris recurrent homogeneous 
Markov process, so, it converges in 
distribution to random vector $\hat X=(\hat Q,\hat U,\hat W)$\,, 
as $t\to\infty$\,, no matter the initial condition,
 see  Dai \cite{dai95}, Asmussen \cite{Asm03}.

Suppose that, furthermore, the following hypothesis holds:
\begin{enumerate}
\item[3.]
  for all
$k\in\mathcal{K}$\,,
$\mathbf{P}(\xi_k=0)=\mathbf{P}(\eta_k=0)=0,$
$\mathbf{E}\exp(\vartheta \xi_k)<\infty$ and 
$\mathbf{E}\exp(\vartheta \eta_k)<\infty$  for some 
$\vartheta>0$, and the c.d.f.'s of 
$\xi_k$ and $\eta_k$ are right--differentiable at $0$ with 
 positive derivatives. 
\end{enumerate}
Let
$\beta_k=\sup\{\vartheta\in\R_+:\,\mathbf{E}\exp(\vartheta
\xi_k)<\infty\}$ and
$\gamma_k=\sup\{\vartheta\in\R_+:\,\mathbf{E}\exp(\vartheta
\eta_k)<\infty\}$\,. 
Let also $\pi(u)=u\ln u-u+1$ if $u>0$\,, $\pi(0)=1$\,,
$\pi(\infty)=\infty$\,,  $0/0=0$\,, and $\infty\cdot 0=0$\,. 
Let $\mathbb{S}_+^{K\times K}$ represent the set of
row--substochastic $K\times K$--matrices.
Given vectors $ a=( a_1,\ldots, a_K)^T\in\R_+^K$
and 
$ d=( d_1,\ldots, d_K)^T\in\R_+^K$,  matrix $  r \in\mathbb{S}_+^{K\times K}$ with
rows $  r _k=(r_{kl}),\,k\in\mathcal{K}$, and $J\subset\mathcal{K}$,
we define 
    \begin{align*}
 \psi^A_k( a_k)&=\sup_{\vartheta<\beta_k}\bl(\vartheta-  a_k\ln \mathbf{E}
\exp(\vartheta \xi_k)\br)\,,\\
\psi^S_k( d_k)&=\sup_{\vartheta<\gamma_k}\bl(\vartheta-  d_k\ln \mathbf{E}
\exp(\vartheta \eta_k)\br)\,,\\
\psi^R_k(  r _{k})&=
\sum_{l=1}^K \pi\Bl(\frac{  r _{kl}}{p_{kl}}\Br)p_{kl}+
\pi\Bl( \frac{1-\sum_{l=1}^K
    r _{kl}}{1-\sum_{l=1}^K p_{kl} }\Br)\bl(1-\sum_{l=1}^K p_{kl}\br)
\end{align*}
and
\begin{equation}
\label{eq:psi_J}
     \psi_J( a, d,  r )= \sum_{k=1}^K\psi_k^A( a_k)+
\sum_{k\in J^c}\psi^S_k( d_k)
+\sum_{k\in J}\psi^S_k( d_k)
\,1_{\{ d_k>\mu_k\}}
+\sum_{k=1}^K  d_k\, \psi^R_k(  r _k),
\end{equation}
where  $J^c=\mathcal{K}\setminus J$\,.
It is readily seen that
 the functions $\psi_k^A( a_k)$\,,
$\psi^S_k( d_k)$ and $\psi^R_k(  r _k)$ are strictly convex,
 the functions
$\psi^S_k( d_k)
\,1_{\{ d_k>\mu_k\}}$ are convex (note that
$\psi^S_k( \mu_k)=0$)  and $\psi_J( a, d,  r )\ge0$\,.
For   $y\in\R^K$, we let
\begin{equation}  \label{eq:Psi}
  L_J(y)=\inf_{\substack{( a, d,  r )\in\R_+^K\times\R_+^K\times
\mathbb{S}_+^{K\times K} 
:\\y= a+(  r ^T-I) d}}\psi_J( a, d,  r )\,.
\end{equation}
  The function $L_J(y)$ is convex by Lemma 4.2 in Puhalskii
  \cite{Puh07}. Besides, the infimum in \eqref{eq:Psi} is attained.
We also note that $L_{\mathcal{K}}(0)=0$ because one can take
$a=\lambda$\,,  $ r =P$ and $d=(I-P^T)^{-1}\lambda$ in \eqref{eq:Psi}.

For  $\,J\subset\mathcal{K}$\,,  
we denote 
$F_J=\{x=(x_1,\ldots,x_K)^T\in \R_+^K:  
x_k=0,k\in J,x_k>0,k\not\in J\}\,$\,.
Let, for
$x\in\R_+^K$ and $y\in\R^K$\,,
\[
  L(x,y)=\sum_{J\subset\mathcal{K}}
1_{ F_J}(x)L_J(y)\,.
\]
As a function of the state
variable $x$\,, the Lagrangian $L(x,y)$ is piecewise constant, being
constant
 on the facets
$F_J$ and having jumps between the facets.
Let \[
  {I}(q)=\int_0^\infty
L(q(t),\dot q(t))\,dt\,,
\]
provided $q=(q(t)\,,t\in\R_+)\in\D(\R_+,\R_+^K)$ is  
absolutely continuous  with
$q(0)=0$ and $I(q)=\infty$\,, otherwise,
 where
 $q(t)=(q_1(t),\ldots,q_K(t))^T$\,. 
  Under the  hypotheses, the  queue length process
$\{(Q(nt)/n\,,t\in\R_+)\,,n\in\N\}$ with $Q(0)=0$ obeys  the 
LDP 
  in $\mathbb D(\R_+,\R^K_+)$ for rate $n$ with
  deviation function $I(q)$\,,
 Puhalskii \cite{Puh07}.
For $x\in\R_+^K$\,, the quasipotential is defined by
\begin{equation}
  \label{eq:4}
  V(x)=\inf_{t\in\R_+}\inf_{\substack{q\in \D(\R_+,\R_+^K):\\
\,q(t)=x}} I(q)\,.
\end{equation}
By Puhalskii \cite{Puh_25b},
the sequence $\{\hat Q/n\,,n\in\N\}$
obeys the LDP in $\R_+^K$ for rate $n$ with the deviation function
 $ V(x)$\,. 
The variational principle ensues,
\[
\lim_{\epsilon\to0}\limsup_{n\to\infty}
\frac{1}{n}\ln\mathbf P\bl(\abs{\frac{\hat Q}{n}-x}<\epsilon\br)=
\lim_{\epsilon\to0}\liminf_{n\to\infty}
\frac{1}{n}\ln\mathbf P\bl(\abs{\frac{\hat Q}{n}-x}<\epsilon\br)=-V(x)\,.
\]
It is  proved in 
Puhalskii \cite{Puh19a} that, for given $x$\,, the infimum in \eqref{eq:4}
 can be taken over some bounded interval $[0,T]$\,, that it is
 attained  and that $ I(q)$ equals $\int_0^TL(q(t),\dot
 q(t))\,dt$\,. 
By convexity, if a trajectory $(q(t),\,a\le t\le
b)$ belongs to $F_J$\,, then
$\int_a^bL(q(t),\dot q(t))\,dt$ is minimised by the straight--line segment
$q(t)=q(a)+(q(b)-q(a))/(b-a)(t-a)$\,. 
It is also
 noteworthy that if  $J'\subset J''$\,,  then $\psi_{J'}(a,d,\varrho)\ge
\psi_{J''}(a,d,\varrho)$\, so that
$L_{J'}(y)\ge
L_{J''}(y)$\,. Hence, if the initial and terminal points of a
piece of an optimal trajectory belong to $F_J$\,, then motion in a straight
line ''costs'' less than motion through adjacent 
regions $F_{J'}$ with $J'\subset J$\,. On the other hand, 
transit from $F_J$ to $F_{J'}$ with $J'\supset J$ takes a positive
length of time and is represented by a straight line segment.
All in all,
an optimal trajectory
to get from the origin to point $x\in\R_+^K$  has to be  piecewise linear.
(We also conjecture that the number of straight line segments that the
optimal trajectory consists of does not
exceed $K$ and that, similarly to a Jackson network,
  zero--component sets $J$ are nonincreasing by
inclusion along the way, with a possible exception of the
terminal point, cf. Puhalskii \cite{Puh19}.)

Let us introduce the Hamiltonian
\begin{equation}
  \label{eq:6}
  H_J(\theta)=\sup_{y\in\R^K}(\theta y-L_J(y))\,,
\end{equation}
where  $\theta=(\theta_1,\ldots,\theta_K)$\,.
($\theta$ being a row vector and $y$ being a column vector, the product
$\theta y$ is well defined.)
We also let $H_0(\theta)=H_{\emptyset}(\theta)$\,.
The function $H_J(\theta)$ is convex, is lower semicontinuous and
$H_J(0)=0$\,. 

In the special case of a Jackson network, i.e. when the interarrival
and service times are exponentially distributed,
it was shown   in Puhalskii  \cite{Puh19} that
\begin{equation}
  \label{eq:15}
  H_J(\theta)=\sum_{k=1}^K (e^{\theta_k}-1)\lambda_k+
\sum_{k\in J}h_k(\theta)^+\mu_k+\sum_{k\in J^c}h_k(\theta)\mu_k\,,
\end{equation}
where the $\lambda_k$ represent  arrival rates, the $\mu_k$
represent  service rates and
\begin{equation}
  \label{eq:4'}
  h_k(\theta)=e^{-\theta_k}\bl(\sum_{l=1}^K(e^{\theta_l}-1)p_{kl}+1\br)-1\,.
\end{equation}
A description of optimal trajectories was also provided.
By Lemma 2 in Puhalskii \cite{Puh19}, an optimal trajectory in 
facet $F_J$
satisfies the equation
\begin{equation}
  \label{eq:7}
  \dot q(t)\in\partial H_J(\theta)\text{ a.e.}
\end{equation}
with constant momentum $\theta$
such that
\begin{equation}
  \label{eq:12}
 H_J(\theta)=0\,,
\end{equation}
 where $\partial$ represents 
subdifferentiation.
 In  each facet an optimal trajectory has to use a certain momentum,
which does not depend on the destination point.
In a generic facet $F_J$\,, the optimal momentum 
$\theta^{(J)}$\,, which is a $K$--vector, is specified by the requirements
that $h_k(\theta^{(J)})=0$\,, if $k\in
J$\,, and $e^{-\theta^{(J)}_k}=    \nu_k/\mu_k$\,, if $k\in J^c$\,,
where $\nu_k$ represents the $k$--th entry of
$\nu=(I-P^T)^{-1}\lambda$\,, with  $\lambda$ representing the vector of
exogeneous arrival rates.
A trajectory to the destination 
is optimal if it uses the optimal momentum
within each facet it traverses and does not visit the same facet more
than once. The trajectory is specified by the time epochs when it
enters  new facets, e.g. if a trajectory follows a coordinate axis
and ''needs'' to enter an adjacent facet, then the only ``control''
 is the turning point on the coordinate axis where that happens.
In order to determine the locations of the turning points,
 we need to start at the
destination and go in reverse time according to the optimal momenta.
This recipe is consistent with the fact that optimal trajectories are
time reversed fluid trajectories of the dual dynamic according to the
 Onsager--Machlup principle.
 Some counter--intuitive properties
were discovered, for instance,  in order to
reach a point on a coordinate axis, it may be advantageous to go down 
another coordinate axis first.

The main result of this paper is the following
generalisation of \eqref{eq:15}.
\begin{theorem}
  \label{the:hamilt}
For $J\subset \mathcal{K}$ and $\theta=(\theta_1,\ldots,\theta_K)$\,,
\begin{equation}
   \label{eq:9}
     H_J(\theta)=\sum_{k=1}^Kf_k(e^{\theta_k}-1)\lambda_k+\sum_{k\in
    J}g_k\bl(
h_k(\theta)^+\br)\mu_k\,
+\sum_{k\in J^c}
g_k\bl(h_k(\theta)\br)\mu_k\,,
\end{equation}
where the functions $f_k(u)$ and $g_k(u)$ are defined for 
$u>-1$ by the respective equations
\begin{equation}
  \label{eq:1a}
  \ln \mathbf E e^{-f_k(e^\vartheta-1)\lambda_k\xi_k}=-\vartheta
\end{equation}
and
\[
  \ln \mathbf E e^{-g_k(e^\vartheta-1)\mu_k\eta_k}=-\vartheta\,,
\]
for  $\vartheta\in\R$\,.
\end{theorem}
Note that $f_k(u)=g_k(u)=u$ in the exponential case.
The proof of the theorem, along with a study of the properties of the
functions $f_k$ and $g_k$\,, is relegated to
Appendix \ref{sec:appendix}.
By \eqref{eq:9}, $H_J(\theta)$ is a locally Lipschitz continuous
function. Hence, \eqref{eq:7}  and \eqref{eq:12} continue to hold for
optimal trajectories.


\section{Two stations}
\label{sec:two-stations}
We assume that $K=2$\,. Then, similarly to Puhalskii \cite{Puh19}, with $\theta=(\theta_1,\theta_2)$\,,
 \begin{align}
  \label{eq:15'}
  h_1(\theta)&=e^{-\theta_1}\bl((e^{\theta_1}-1)p_{11}+(e^{\theta_2}-1)p_{12}+1\br)-1\,,\\
\notag
  h_2(\theta)&=e^{-\theta_2}\bl((e^{\theta_1}-1)p_{21}+(e^{\theta_2}-1)p_{22}+1\br)-1\,\,,\\
\notag  H_1(\theta)&=f_1(e^{\theta_1}-1)\lambda_1+f_2(e^{\theta_2}-1)\lambda_2+
g_1(h_1(\theta)^+)\mu_1+g_2(h_2(\theta))\mu_2,\\
\label{eq:15b}
  H_2(\theta)&=f_1(e^{\theta_1}-1)\lambda_1+f_2(e^{\theta_2}-1)\lambda_2+
g_1(h_1(\theta))\mu_1+g_2(h_2(\theta)^+)\mu_2,
\\\notag
 H_0(\theta)&=f_1(e^{\theta_1}-1)\lambda_1+f_2(e^{\theta_2}-1)\lambda_2+
g_1(h_1(\theta))\mu_1
+g_2(h_2(\theta))\mu_2\,.
\end{align}

Some  thought shows that an optimal trajectory to get from the
origin to
$x\in\R_+^2$ must either go down one of the coordinate axes first and
then transit to $x$ or go to $x$ in a straight line from the outset. 
(For some insight, see
Appendix \ref{sec:optim-two-stra}.)
 Suppose that an optimal piecewise linear
trajectory $q(t)=(q_1(t),q_2(t))^T$ gets from the origin to $x=(x_1,x_2)^T$ 
with $x_2>0$ by following first the
$x_1$--axis and going next in
a straight line to $x$\,. Let $(s^{(1)}_1,0)$  represent
 the point on the $x_1$--axis where $q(t)$
leaves the $x_1$--axis. Let $\theta^{(1)}=
(\theta_1^{(1)},\theta_2^{(1)})$ and 
$\theta^{(1,2)}=(\theta_1^{(1,2)},\theta_2^{(1,2)})$ 
represent the momenta of going down the $x_1$--axis 
and traversing the interior of the $x_1x_2$--quarter plane,
respectively.
\begin{theorem}
  \label{2D}
 Suppose that    $s^{(1)}_1>0$\,.
The momentum $\theta^{(1)}$ is specified by the
equations 
$  h_2(\theta^{(1)})=0$ and
$  H_0(\theta^{(1)})=0$\,.
 The
 momentum $\theta^{(1,2)}$ is specified by the
equations 
\begin{equation}
  \label{eq:31'}
  \theta^{(1,2)}_1=\theta_1^{(1)}
\end{equation}
and $H_0(\theta^{(1,2)})=0\,.$
In addition, \begin{equation}
  \label{eq:29'}
    s^{(1)}_1=x_1-\frac{\partial_1H_0(\theta^{(1,2)})}{\partial_2H_0(\theta^{(1,2)})}\,x_2\,.
\end{equation}
\end{theorem}

\begin{proof}
  The motion down the $x_1$--axis occurs
with Hamiltonian $H_{2}$ under the constraint that $q_2(t)=0$\,.
We have that, a.e.,
 $ (\dot q_1(t),0)^T\in\partial H_{2}(\theta^{(1)})$
 so that
 $0\in\partial_2 H_{2} (\theta^{(1)})$\,.
Since  by \eqref{eq:15'}, $\partial_2
h_1(\theta ) = e^{-\theta_1+\theta_2}p_{12}\ge0$\,,  
\eqref{eq:15b} implies that in order for the latter
inclusion to hold
we must have that $h_2( \theta^{(1)})\ge0$\,. Hence, 
$  H_0(\theta^{(1)})=H_2(\theta^{(1)})=0\,.$
The geometric argument of Section 5  in Puhalskii
\cite{Puh19} shows that if $h_2(\theta^{(1)})>0$\,, then going 
from the origin to
 $x$ in a straight line costs less than first
going to  $(s_1^{(1)},0)^T$ with $\theta^{(1)}$\,.
Therefore, $h_2(\theta^{(1)})=0$\,.

Let $ p=(p_1,p_2)^T$ represent the  displacement vector in the $x_1x_2$--plane:
$ p= x- s^{(1)}$\,, where $s^{(1)}=(s^{(1)}_1,0)^T$\,.
The motion from $s^{(1)}$ to $x$ occurs  with
momentum $\theta^{(1,2)} =(\theta^{(1,2)}_1,\theta^{(1,2)}_2)$ 
such that $H_0(\theta^{(1,2)})=0$ and the outward normal $\nabla H_0(\theta^{(1,2)})$ to the curve
$H_0(\cdot)=0$ at $\theta^{(1,2)}$ has the same slope as $p$\,, see
\eqref{eq:7}.  Since $p_2=x_2>0$\,, $\partial_2
H_0(\theta^{(1,2)})>0\,.$
Therefore, for all $\theta=(\theta_1,\theta_2)$ 
 in a suitable neighbourhood of $\theta^{(1,2)}$ such that 
$H_0(\theta)=0$\,,
$\theta_2$ can be expressed
   as a
strictly concave twice continuously differentiable function of $\theta_1$ and 
\begin{equation}
  \label{eq:43}
  \theta_2'=-\frac{\partial_1H_0(\theta)}{\partial_2 H_0(\theta)}\,,
\end{equation}
 with
prime referring to  $\theta_1$--differentiation.
The slope of $p$ is $\partial_2H_0(\theta)/\partial_1 H_0(\theta)=-1/\theta_2'$ so that
\begin{equation}
  \label{eq:51'}
   p_1=-\theta_2' p_2=-\theta_2' x_2\,.
\end{equation}

We have that
\begin{equation}
  \label{eq:53}
    s^{(1)}_1=x_1-p_1=x_1-\frac{\partial_1H_0(\theta^{(1,2)})}{\partial_2H_0(\theta^{(1,2)})}\,x_2\,,
\end{equation}
which is \eqref{eq:29'}.

Let $\theta(x)=(\theta_1(x),\theta_2(x))$ denote
the momentum of getting to $x$ in a straight line. Equivalently, it represents
 the point on the curve $H_0(\cdot)=0$ such that
the outward normal at that point is colinear to $x$ \,, i.e.
\[
  \frac{\partial_1H_0(\theta(x))}{\partial_2H_0(\theta(x))}=\frac{x_1}{x_2}\,. 
\]
Thus,
\[
  \theta_2'(\theta_1(x))x_2=-x_1.
\]
The existence of $\theta(x)$ holds because the Gauss map on a compact
oriented manifold is surjective, see, e.g. Thorpe \cite{tho79}. In
addition, uniqueness holds because the function $\theta_2$ is strictly
concave so that its derivative is injective.

As $s^{(1)}_1>0$ and by  \eqref{eq:43} and \eqref{eq:53},
\begin{equation}
  \label{eq:63}
  s^{(1)}_1
=x_1+\theta_2'( \theta^{(1,2)}_1)\,x_2\,,
\end{equation}
we have that
  $x_1>-\theta_2'(\theta^{(1,2)}_1)x_2$ so that,
due to $-\theta_2'(\cdot)$ being strictly increasing,
\begin{equation}
  \label{eq:30'}
  \theta^{(1,2)}_1< \theta_1(x)\,.
\end{equation}

As in  (2.19) and (2.20) in Puhalskii \cite{Puh19},
the total cost of getting to $x$ provided momentum $\theta$ is used in the
$x_1x_2$--plane  is
\begin{equation}
  \label{eq:40}
    \theta_1^{(1)}s^{(1)}_1+\theta p
=\theta_1^{(1)}x_{1}-\theta_1^{(1)} p_{1}
+\theta_1 p_1+\theta_2 p_2\,.
\end{equation} 
On taking into account \eqref{eq:51'}, we have 
 that 
 $\theta_1^{(1,2)}$  minimises 
$-\theta_1^{(1)} p_{1}
+\theta_1 p_1+\theta_2 p_2=\bl( -(\theta_1
-\theta_1^{(1)})\theta_2'+\theta_2\br)p_2
=\bl( -(\theta_1
-\theta_1^{(1)})\theta_2'+\theta_2\br)x_2$ as a function of $\theta_1$
over $\theta_1< \theta_1(x)$\,.
The $\theta_1$--derivative of this expression is
$-(\theta_1-\theta^{(1)}_1)\theta_2''x_2$\,.
By strict
concavity, $\theta_2''<0$\,.
Since  $\theta_1^{(1,2)}$ delivers a minimum and 
\eqref{eq:30'} holds,  the derivative at that
point is zero so that $\theta_1^{(1,2)}=\theta_1^{(1)}$\,.

\end{proof}
\begin{figure}[h!]
  \centering
  \includegraphics[scale=.7]{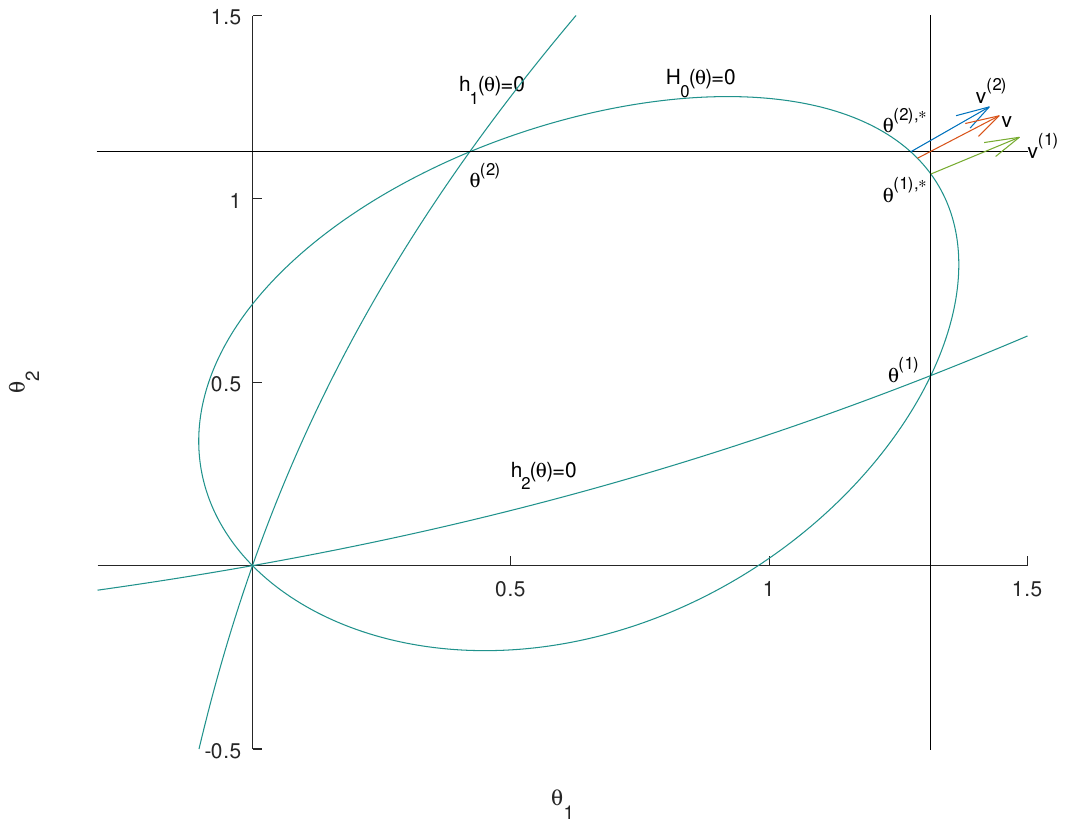}
  \caption{Intersection outside the domain}
  \label{fig:1}
\end{figure}
\begin{remark}
  The equality \eqref{eq:31'} constitutes a conservation law:
 if an
optimal trajectory spends time in
 the facet $F_2=\{x_1>0\,, x_2=0\}$\,, i.e.
    $s^{(1)}_1>0$\,,
then the $\theta_1$--component of its continuation to the facet
$F_0=\{x_1>0\,, x_2>0\}$ is the same as the
$\theta_1$--component
in $F_2$\,.
\end{remark}
\begin{remark}
  By \eqref{eq:40} and conservation of momenta,
under the hypotheses the cost of getting to $x$ is
$V(x)=\theta^{(1,2)}x$\,. 

\end{remark}
\begin{remark}
  The proof implies that, generally, 
$\theta_1^{(1,2)}=\theta^{(1)}_1\wedge  \theta_1(x)\,,$
for all $s^{(1)}_1\ge0$\,.
\end{remark}
\begin{remark}
If $\theta=(\theta_1,\theta_2)$ is a generic point on the curve $H_0(\cdot)=0$\,,
then
the slope of $\nabla H_0(\theta)$
is greater than the slope of $x$ if and only if $x_1>-\theta_2'(\theta_1)x_2$\,.
\end{remark}
\begin{remark}
 We have that $s^{(1)}_1>0$ if only if
 $\theta^{(1)}_1< \theta_1(x)$
 (necessity follows from the proof, sufficiency uses \eqref
 {eq:63}). 
Another equivalent condition is that the slope of $x$ is less than
the slope of $\nabla
H_0(\theta^{(1),\ast})$\,,
where $\theta^{(1),\ast}=(\theta_1^{(1)},\theta_2(\theta_1^{(1)}))$\,.
\end{remark}
\begin{remark}
 The  equation in \eqref{eq:51'} can also be derived from the requirement that
$\theta$ deliver $    \sup_{\theta:\, H_0(\theta)=0}\theta p$ so that
the full $\theta_1$--derivative of $\theta p$ on the set
$H_0(\theta)=0$ is zero. A similar
line of reasoning will be pursued in the next section
 when analysing a 3--station network.
Besides,  the equation in question is equivalent to the Hamilton equation
$   p_1=\partial_1H_0(\theta)/
\partial_2H_0(\theta)\,x_2\,.
$\end{remark}

\begin{figure}[h!]
  \centering
  \includegraphics[scale=.7]{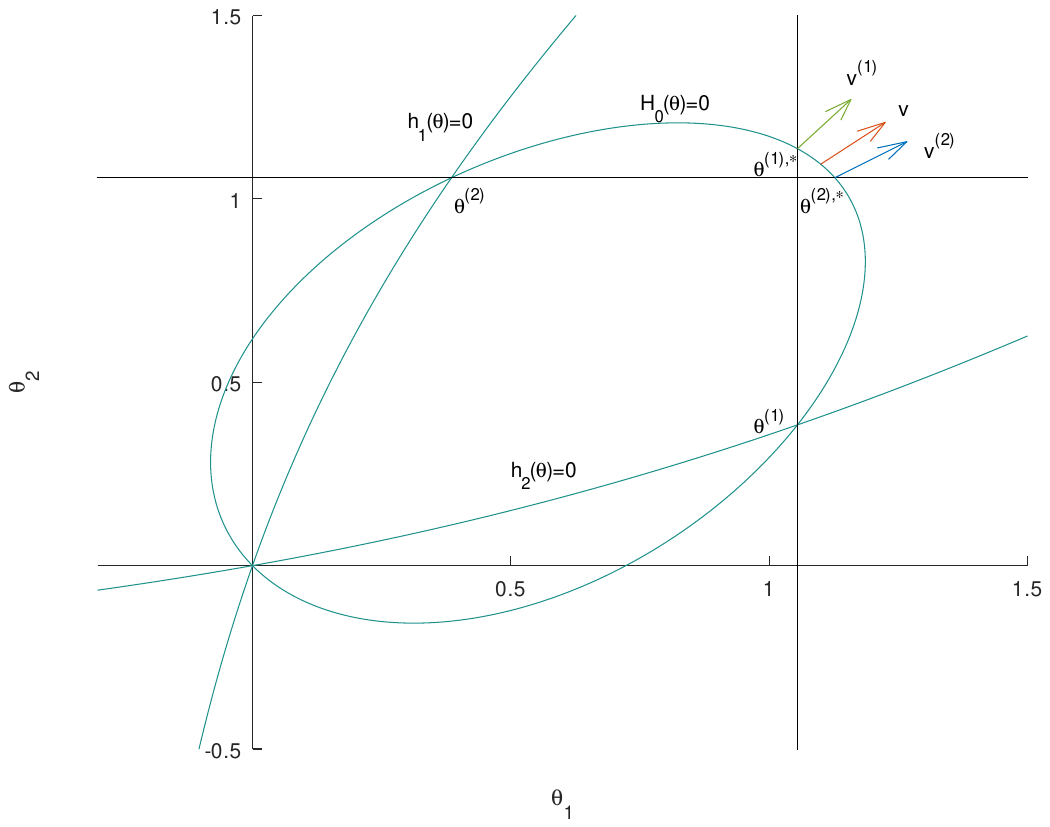}
  \caption{Intersection inside the domain}
  \label{fig:2}
\end{figure}

We illustrate the above analytical considerations with graphs which
bring into play the other coordinate axis. The exposition is in a
similar vein as the contents of Section 5 in Puhalskii \cite{Puh19}.
As before, let $\theta^{(1)}$ represent the  intersection  of the curves
$H_0(\theta)=0$ and $h_2(\theta)=0$ and let 
$\theta^{(2)}$ represent the intersection  of the curves $H_0(\theta)=0$
and $h_1(\theta)=0$\,. Let $\theta^{(1),\ast}$\,, respectively, $\theta^{(2),\ast}$
represent the point $\theta=(\theta_1,\theta_2)$ on the curve
 $H_0(\theta)=0$ such that $\theta_1=\theta^{(1)}_1$\,, respectively,
such that $\theta_2=\theta^{(2)}_2$\,,  see
Fig.\ref{fig:1}. 
Let $v^{(1)}$ and $v^{(2)}$ represent 
unit outward normals to $H_0(\theta)=0$ at 
$\theta^{(1),\ast}$ and $\theta^{(2),\ast}$\,, respectively.
Suppose that  the slope  of vector 
$x$ 
is less than  the slope of $v^{(1)}$\,(which is 
$\partial_2H_0(\theta^{1,\ast})
/\partial_1H_0(\theta^{1,\ast})=-1/\theta_2'(\theta^{(1)}_1)$)\,.
Then, $\theta^{(1)}_1< \theta_1(x)$\,.
 When going down the $x_1$--axis,
it is optimal to use momentum
 $\theta^{(1)}$ until the direction to $x$ is given by
$v^{(1)}$\,. 
From that point on,
it is optimal  to go toward $x$ in a straight line 
(with $\theta^{(1),\ast}$)\,. The total  cost is
$\theta^{(1),\ast}\, x$\,.
Indeed, the cost of  displacement by  $s^{(1)}_1$ along
the horizontal axis with $\theta^{(1)}$
 followed by going straight to $x$ with
$\theta^{(1),\ast}$ is 
$\theta^{(1)}_1s^{(1)}_1+\theta^{(1),\ast}(x-(s^{(1)}_1,0)^T)=\theta^{(1),\ast}x\,$.   Following that trajectory is superior to going to
$x$ in a straight line. On the other hand, if 
the slope of  $x$ is greater
than  the slope of $v^{(1)}$\,,
 then going in a straight line
  to $x$ (with $ \theta(x)$) is superior to
going horizontally initially and $\theta^{(1)}_1>\theta_1(x)$. 
Switching the roles of $x_1$ and $x_2$ implies that if the slope of $x$ 
is greater than  the slope of $v^{(2)}$\,, then,
going up the vertical axis with $\theta^{(2)}$
  until the direction to $x$ is given by
$v^{(2)}$ and going in a straight line to $x$  with
$\theta^{(2),\ast}$ afterwards
 is superior to going to $x$ directly and is optimal
among all  paths that start off vertically.
 The  cost is
$\theta^{(2),\ast} x$\,. If the slope of vector $x$ is less
than the slope of $v^{(2)}$\,, then going straight to $x$ is superior to
going vertically initially.

Since the slope of $x$ is less than  the
  slope of $v^{(i)}$
if and only if $\theta(x)$ is clockwise w.r.t.
$\theta^{(i),\ast}$\,, where $i\in\{1,2\}$\,, we obtain that
\begin{figure}[h!]
  \centering
  \includegraphics[scale=.7]{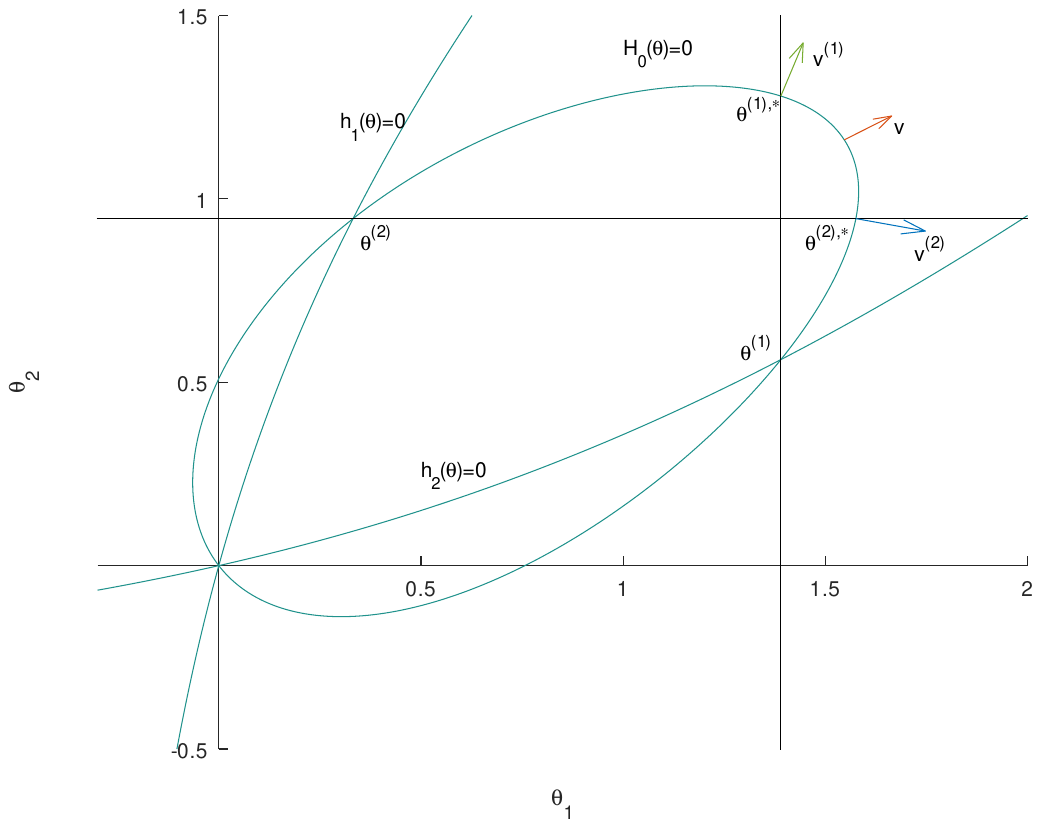}
  \caption{Intersection inside the domain}
  \label{fig:3}
\end{figure}
\begin{figure}[h!]
  \centering
  \includegraphics[scale=.7]{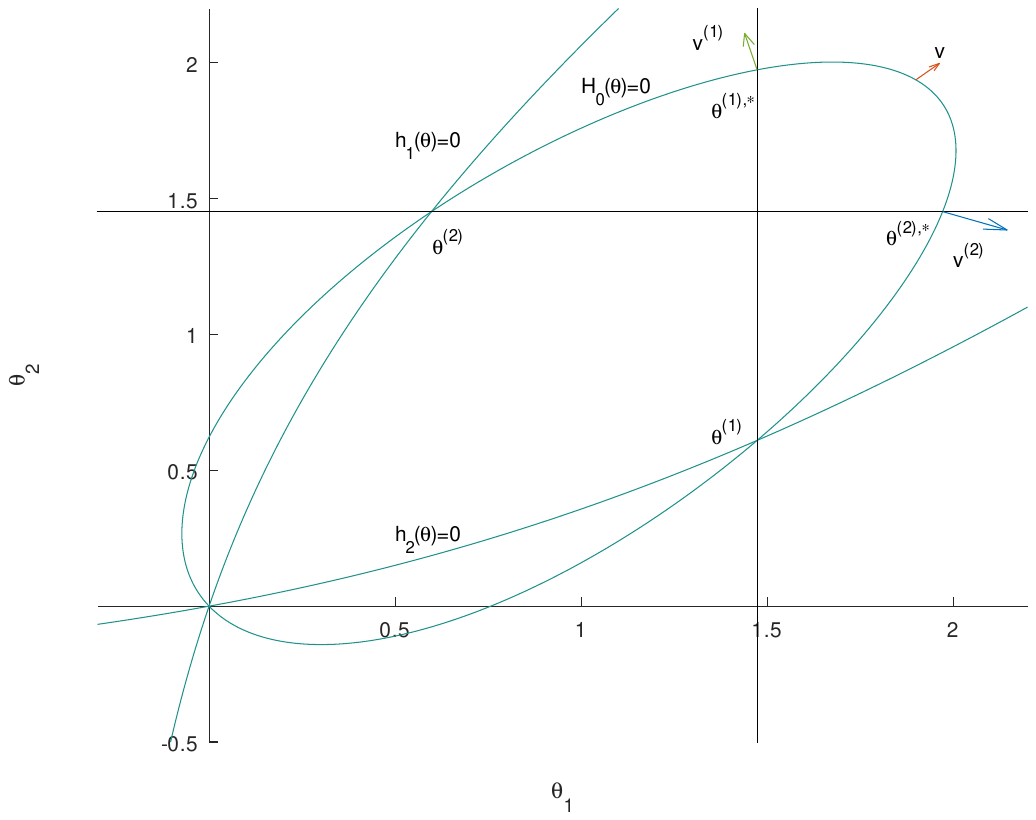}
  \caption{Intersection inside the domain}
  \label{fig:4}
\end{figure}
\[
  V(x)=
  \begin{cases}
\theta^{(1),\ast}\, x\,, &\text{ if 
$\theta(x)$ is clockwise w.r.t.
$\theta^{(1),\ast}$}\\
\theta^{(2),\ast}\, x\,, &\text{
if 
$\theta(x)$ is anticlockwise w.r.t.
$\theta^{(2),\ast}$ }\\
\theta(x)\,x\,, &\text{ if $\theta(x)$ is on the arc connecting 
$\theta^{(2),\ast}$ and $\theta^{(1),\ast}$
}
      \end{cases}
\]
 In Fig.\ref{fig:1}
the straight lines $\theta_1=\theta_1^{(1)}$ and $\theta_2=\theta_2^{(2)}$
intersect outside  the domain
$\{\theta:\,H_0(\theta)\le0\}$ so that  $\theta^{(2),\ast}$ is
anticlockwise w.r.t. $\theta^{(1),\ast}$ 
and the slope of $v^{(2)}$
 is   greater than the slope  of $v^{(1)}$\,.
 (Recall that for a Jackson network
the intersection point is on the line $H_0(\theta)=0$\,, Puhalskii
\cite{Puh19}.)

Suppose now that 
the straight lines $\theta_1=\theta_1^{(1)}$ and $\theta_2=\theta_2^{(2)}$
intersect inside  the domain $\{\theta:\,H_0(\theta)\le0\}$ so that
 $\theta^{(2),\ast}$ is
clockwise w.r.t.  $\theta^{(1),\ast}$ 
and  the slope of $v^{(2)}$
 is   less than the slope  of $v^{(1)}$\,, see Fig.\ref{fig:2}.
If $\theta(x)$ is clockwise 
w.r.t. $\theta^{(1),\ast}$ (i.e. the slope of $x$ is less than 
 the slope of $v^{(1)}$), then 
 going horizontally first (with $\theta^{(1)}$) and following the direction of $v^{(1)}$
next is superior to going to $x$ in a straight line. For points $x$ such that
 $\theta(x)$ is  anticlockwise  w.r.t.
$\theta^{(1),\ast}$ going to $x$ in a straight line (with $\theta(x)$)
is 
superior to 
going horizontally initially. Similarly, going
vertically first (with $\theta^{(2)}$) and then following the
direction of $v^{(2)}$ (with $\theta^{(2),\ast}$) is
superior to going in a straight line to $x$ provided $\theta(x)$ is
anticlockwise w.r.t.
 $\theta^{(2),\ast}$ and is inferior to going in a straight line to
$x$ provided $\theta(x)$ is 
clockwise w.r.t. $\theta^{(2),\ast}$\,.  

Thus, if $\theta(x)$ is anticlockwise w.r.t.
$\theta^{(1),\ast}$ (so, anticlockwise w.r.t. $\theta^{(2),\ast}$), it is optimal to go vertically until the
slope to $x$ is equal to the slope of $v^{(2)}$ and then go in a straight
line to $x$\,, the cost being $\theta^{(2),\ast}x$\,.
If  $\theta(x)$ is clockwise w.r.t.
$\theta^{(2),\ast}$ (so, clockwise w.r.t. $\theta^{(1),\ast}$), then it is optimal to go horizontally until the
slope to $x$ is equal to the slope of $v^{(1)}$ and then go straight
to $x$\,, the cost being $\theta^{(1),\ast}x$\,.

Suppose that   $\theta(x)$  is on the arc connecting
$\theta^{(1),\ast}$ and $\theta^{(2),\ast}$
(i.e. the slope of $x$ is  between the slopes
of $v^{(1)}$ and $v^{(2)}$)\,. Then,
there are two competing trajectories. 
The first one is
down the horizontal axis and then colinearly to $v^{(1)}$\,, the cost being 
$\theta^{(1),\ast} x$\,.
The other one  is 
up the vertical axis and then colinearly to $v^{(2)}$\,,
 the cost being 
$\theta^{(2),\ast} x$\,.
The
horizontal path should be used if $(\theta^{(1),\ast}-\theta^{(2),\ast}) x<0$
and the vertical path should be used if
$(\theta^{(1),\ast}-\theta^{(2),\ast}) x>0$\,.
Let $v_0$ represent the unit outward normal to the curve $H_0(\theta)=0$
 that is orthogonal to
$\theta^{(1),\ast}-\theta^{(2),\ast}$ and let $\theta(v_0)$ denote the
associated point. (The slope of $v_0$ is  between the slopes of
$v^{(1)}$ and $v^{(2)}$\,.) If $\theta(x)$ is clockwise w.r.t.
 $\theta(v_0)$\,, then the horizontal path should
be used. If $\theta(x)$ is anticlockwise w.r.t.
$\theta(v_0)$\,, then 
the vertical path should be used.
The analysis can be summarised by stating that for the configuration in Fig.2,
\[
    V(x)=
  \begin{cases}
\theta^{(1),\ast}\, x\,, &\text{ if 
$\theta(x)$ is clockwise w.r.t. $\theta(v_0)$
}
\\
\theta^{(2),\ast}\, x\,, &\text{ if $\theta(x)$ is 
  anticlockwise  w.r.t. $\theta(v_0)$
}
      \end{cases}
\]
Some other  cases are depicted in Fig.3 and Fig.4.
It is noteworthy that going in a straight line to $x$  can be optimal only
if $\theta^{2,\ast}$ is anticlockwise w.r.t. $\theta^{1,\ast}$\,, as in
Fig. 1.

\section{Three stations}
\label{sec:three-stations}
We now analyse a three station network. 
Suppose that an optimal piecewise linear
trajectory $q(t)=(q_1(t),q_2(t),q_3(t))^T$ gets from the origin to 
$x=(x_1,x_2,x_3)^T$ 
with $x_3>0$ by going down the
$x_1$--axis,  traversing the  $x_1x_2$--plane,
and going   to $x$\,. Let $s^{(1)}=(s^{(1)}_1,0,0)$  represent
 the point  where $q(t)$
leaves the $x_1$--axis and let $s^{(1,2)}=(s^{(1,2)}_1,s^{(1,2)}_2,0)$  represent
 the point  where $q(t)$
leaves the $x_1x_2$--plane.  Let $\theta^{(1)}=
(\theta_1^{(1)},\theta_2^{(1)},\theta_3^{(1)})$\,,
$\theta^{(1,2)}=(\theta_1^{(1,2)},\theta_2^{(1,2)},\theta_3^{(1,2)})$
and $\theta^{(1,2,3)}=(\theta_1^{(1,2,3)},\theta_2^{(1,2,3)},
\theta_3^{(1,2,3)})$ represent the respective
 momenta of going down the $x_1$--axis,
 traversing  the $x_1x_2$--plane, and 
 going from the $x_1x_2$--plane to $x$\,.
\begin{theorem}
  \label{3D}
Suppose that  $s^{(1)}_1>0$\,, $s^{(1,2)}_1>0$ and
$s^{(1,2)}_2>0$\,.
The momentum $\theta^{(1)}$ is specified by the
equations 
$  h_2(\theta^{(1)})=0$\,,
$  h_3(\theta^{(1)})=0$ and
$  H_0(\theta^{(1)})=0$\,,
 the momentum $\theta^{(1,2)}$ is specified by the
equations $\theta^{(1,2)}_1=\theta_1^{(1)}$\,, $h_3(\theta^{(1,2)})=0$\,,
and $H_0(\theta^{(1,2)})=0\,,$ and
the momentum $\theta^{(1,2,3)}$   is
specified by the equations $\theta_i^{(1,2,3)}=\theta_i^{(1,2)}$\,,
for $i\in\{1,2\}$\,, and $H_0(\theta^{(1,2,3)})=0$\,.
The 
 turnoff points $(s^{(1)}_1,0,0)$ and $(s^{(1,2)}_1,s^{(1,2)}_2,0)$
 are
  as follows,
\[
  s^{(1)}_1=s^{(1,2)}_1-
\frac{\partial_1H_0(\theta^{(1,2)})}{\partial_2H_0(\theta^{(1,2)})}\,
s^{(1,2)}_2
\]
and
\begin{equation}
  \label{eq:29}
    s^{(1,2)}_i=x_i-
\frac{\partial_iH_0(\theta^{(1,2,3)})}{\partial_3H_0(\theta^{(1,2,3)})}\,x_3\,,
\end{equation}
for $i\in\{1,2\}$\,.
\end{theorem}
\begin{proof}
The proof employs the methods used for the 2--station network.
  The motion down the $x_1$--axis occurs
with Hamiltonian $H_{\{2,3\}}$ under the constraints that
$q_2(t)=q_3(t)=0$ and is governed by the equation 
$(\dot q_1(t),0,0)^T\in\partial H_{\{2,3\}}(\theta^{(1)})$\,.
The  discusssion of optimal motion down an axis in a general GJN to be
 found in 
Appendix \ref{sec:optimal-paths} implies the equations in the statement.
The motion in  the $x_1x_2$--plane occurs
with Hamiltonian $H_3$ under the constraint that $q_3(t)=0$ and is governed by the equation $(\dot q_1(t),\dot q_2(t),0)^T\in\partial H_{3}(\theta^{(1,2)})$\,.
A similar  argument as in Section \ref{sec:two-stations}
 shows that 
\begin{equation}
  \label{eq:35}
  h_3(\theta^{(1,2)})=0\,.
\end{equation}
Then,
\begin{equation}
  \label{eq:36}
  H_0(\theta^{(1,2)})= H_3(\theta^{(1,2)})=0\,.
\end{equation}
 Let $ p^{(1,2)}$ represent the  displacement vector in the
$x_1x_2$--plane, i.e.
$ p^{(1,2)}= s^{(1,2)}-s^{(1)}$\,.
In view of \eqref{eq:36},
arguing as in (2.19)  in Puhalskii \cite{Puh19}
and noting that the third entry of $p^{(1,2)}$ equals zero obtains that
$\theta^{(1,2)}$ delivers 
\begin{multline*}
    \sup_\theta(\theta  p^{(1,2)}-H_3(\theta))
=\sup_{\theta:\,h_3(\theta)=0}(\theta  p^{(1,2)}-H_3(\theta))
=\sup_{\theta:\,h_3(\theta)=0}(\theta  p^{(1,2)}-H_0(\theta))\\=
\sup_\theta(\theta  p^{(1,2)}-H_0(\theta))=
    \sup_{\tilde\theta}(\tilde\theta \tilde p^{(1,2)}-
\tilde H_0(\tilde\theta))\,,
\end{multline*}
where $\tilde\theta=(\theta_1,\theta_2)$\,, 
   $\tilde p^{(1,2)}=( p_1^{(1,2)}, p_2^{(1,2)})$\,
and $\tilde H_0(\tilde\theta)=\inf_{\theta_3}H_0(\theta)$\,.
The function $ H_0$ being strictly convex, 
the infimum in the definition of $\tilde H_0$ is attained at
a unique $\theta_3$\,, which we denote by
$\theta_3(\tilde\theta)$\,. The equation for the latter function is as follows,
\[
  \partial_3H_0(\tilde\theta,\theta_3(\tilde\theta))=0\,.
\]
By strict convexity,
$\partial_3^2H_0(\tilde\theta,\theta_3(\tilde\theta))>0$ so that, by the
implicit function theorem, $\theta_3(\tilde\theta)$ is a twice
differentiable function. It follows that $\tilde H_0(\tilde\theta)=
H_0(\tilde\theta,\theta_3(\tilde\theta))$
is a twice differentiable function.
Taking into account that $s_2^{(1,2)}>0$ and reasoning 
 as in the proof of Theorem \ref{2D}, we have that in a
neighbourhood of $(\theta_1^{(1,2)},\theta_2^{(1,2)})$ we can write $\theta_2$ as a
twice differentiable function of  $\theta_1$ provided $\tilde H_0(\tilde\theta)=0$\,. 
It follows
that $(\theta_2,\theta_3)$ is a twice differentiable function of
$\theta_1$ locally on the set $\{\theta:\,H_0(\theta)=0\}$\,, i.e.
$(\theta_2,\theta_3)=f^{(1,2)}(\theta_1)$
 with
 $f^{(1,2)}$ being twice
differentiable and strictly concave.

Using 
prime in order to denote  $\theta_1$--differentiation, we have that if
$\sup_{\theta}(\theta  p^{(1,2)}-H_0(\theta))=
\sup_{\theta}\theta  p^{(1,2)}$
 is attained at $\theta$\,, then
$0=d/d\theta_1\,(\theta  p^{(1,2)})=   p^{(1,2)}_1+\theta_2'
p^{(1,2)}_2\,,$
where  $ p^{(1,2)}=( p^{(1,2)}_1, p^{(1,2)}_2,0)^T$\,.
Hence, similarly to
\eqref{eq:51'}, 
\[
   p^{(1,2)}_1=-\theta_2' p^{(1,2)}_2=-\theta_2' s^{(1,2)}_2\,.
\]
In analogy with the notation of Section \ref{sec:two-stations}, 
let $\theta(y)$\,, for $y\in\R_+^3$\,, denote the momentum of getting to $y$ in a straight
line, i.e. $H_0(\theta(y))=0$ and
$ \nabla H_0(\theta(y))$ is colinear to $y$\,.
We  have  analogues of \eqref{eq:63} and \eqref{eq:30'} to the effect that
\[
  s^{(1)}_1=s^{(1,2)}_1+\theta_2'( \theta_1^{(1,2)})\,s^{(1,2)}_2
\]
and
\begin{equation}
  \label{eq:30''}
  \theta^{(1,2)}_1< \theta_1(s^{(1,2)})\,,
\end{equation}
with $\theta_1(s^{(1,2)}) $ being specified by the equation
\begin{equation}
  \label{eq:22'}
  -\theta_2'(\theta_1(s^{(1,2)}))
s^{(1,2)}_2=s^{(1,2)}_1\,.
\end{equation}

The cost of getting to  $s^{(1,2)}$ with momentum $\theta=(\theta_1,\theta_2)$ is
\[
    \theta_1^{(1)}s_1^{(1)}+\theta p^{(1,2)}
=\theta_1^{(1)}s^{(1,2)}_{1}-\theta_1^{(1)} p^{(1,2)}_{1}
+\theta_1 p^{(1,2)}_1+\theta_2 p^{(1,2)}_2
=\theta_1^{(1)}s^{(1,2)}_{1}-(\theta_1-\theta_1^{(1)}) \theta_2'
s^{(1,2)}_2
+\theta_2 s^{(1,2)}_2\,.
\]
The $\theta_1$--derivative is $-(\theta_1-\theta_1^{(1)}) \theta_2''s^{(1,2)}_2$\,.
Since $\theta^{(1,2)}$ is optimal and
\eqref{eq:30''} holds,  the derivative at $\theta_1^{(1,2)}$ must
 equal zero so that $\theta^{(1,2)}_1=\theta_1^{(1)}$\,.
The other two entries of $\theta^{(1,2)}$ are  determined by 
equations  \eqref{eq:35} and \eqref{eq:36}.

Let $p^{(1,2,3)}=x-s^{(1,2)}$\,.
We note that $p_3^{(1,2,3)}=x_3$\,.
Since $H_0(\theta^{(1,2,3)})=0$\,, it may be assumed that locally
 $\theta_3^{(1,2,3)} 
=f^{(1,2,3)}(\theta_1^{(1,2,3)} ,\theta_2^{(1,2,3)} )\,,$ with 
$f^{(1,2,3)}$ being twice continuously differentiable and strictly concave.
In analogy with the developments in Section \ref{sec:two-stations},
 we now minimise $\theta p^{(1,2,3)}$ by setting the respective partial 
$\theta_i$--derivatives of $\theta p^{(1,2,3)}$ equal to zero:
\begin{equation}
  \label{eq:54}
  p_i^{(1,2,3)}+\partial_i \theta_3x_3 
=0\,,\;
\end{equation}
where $\theta_3(\theta_1,\theta_2)=f^{(1,2,3)}(\theta_1 ,\theta_2 )$
and $i\in\{1,2\}$\,.
In analogy with the developments above, 
\[
s^{(1,2)}_i=x_i-p^{(1,2,3)}_i=x_i+\partial_i \theta_3x_3\,,
\]
so that  $s^{(1,2)}_i>0$ if and only if
\begin{equation}
  \label{eq:52}
  \theta_i^{(1,2,3)}< \theta_i(x) \,,
\end{equation}
where 
\[
  -\partial_i \theta_3(\theta_1(x),\theta_2(x))
=  \frac{x_i}{x_3}\,.
\]
(Generally,  $\theta_i^{(1,2,3)}\le \theta_i(x)$\,.)

The total cost of getting to $x$ is
\begin{multline}
  \label{eq:37}
  \theta_1^{(1)}s^{(1)}_1+\theta^{(1,2)}p^{(1,2)}+\theta^{(1,2,3)}p^{(1,2,3)}=
  \theta_1^{(1)}(x_1-p_1^{(1,2,3)}-p_1^{(1,2)})+\theta^{(1,2)}p^{(1,2)}+\theta^{(1,2,3)}p^{(1,2,3)}\\=
\theta_1^{(1)}x_1+(\theta^{(1,2)}_1-\theta_1^{(1)})p_1^{(1,2)}+
\theta^{(1,2)}_2p^{(1,2)}_2+
(\theta^{(1,2,3)}_1-\theta^{(1)}_1)p_1^{(1,2,3)}
+\theta_2^{(1,2,3)}p_2^{(1,2,3)}+\theta_3^{(1,2,3)}p_3^{(1,2,3)}\\
=\theta_1^{(1)}x_1+(\theta^{(1,2)}_1-\theta_1^{(1)})p_1^{(1,2)}+
\theta^{(1,2)}_2x_2-\theta^{(1,2)}_2p^{(1,2,3)}_2+
(\theta^{(1,2,3)}_1-\theta^{(1)}_1)p_1^{(1,2,3)}
+\theta_2^{(1,2,3)}p_2^{(1,2,3)}+\theta_3^{(1,2,3)}p_3^{(1,2,3)}\\
=\theta_1^{(1)}x_1+(\theta^{(1,2)}_1-\theta_1^{(1)})p_1^{(1,2)}+
\theta^{(1,2)}_2x_2+(\theta_2^{(1,2,3)}-\theta^{(1,2)}_2)p^{(1,2,3)}_2+
(\theta^{(1,2,3)}_1-\theta^{(1)}_1)p_1^{(1,2,3)}
+\theta_3^{(1,2,3)}x_3\,. 
\end{multline}

Owing to  \eqref{eq:31'}, \eqref{eq:54} and \eqref{eq:37},  $(\theta_1^{(1,2,3)},
\theta_2^{(1,2,3)})$ minimises
$
-(\theta_1-\theta^{(1)}_1)\partial_1\theta_3
-(\theta_2      -\theta^{(1,2)}_2)\partial_2\theta_3+\theta_3
$
over $(\theta_1,\theta_2)$ such that $\theta_i<
\theta_i(x)$\,.
The respective 
partial $\theta_1$-- and $\theta_2$--derivatives are as follows,
\[
-(\theta_1-\theta_1^{(1)})\partial_1^2\theta_3
-  (\theta_2-\theta^{(1,2)}_2)\partial_{12}\theta_3
\]
and
\[
-(\theta_1-\theta^{(1)}_1)\partial_{12}\theta_3
- (\theta_2-\theta^{(1,2)}_2)\partial^2_2\theta_3
\,.
\]
Thus, the gradient of the cost is
\[
  -\left(
    \begin{array}{cc}
    \partial_1^2\theta_3&\partial_{12}\theta_3\\
\partial_{12}\theta_3&\partial_2^2\theta_3
    \end{array}
\right)\left(
  \begin{array}[c]{c}
    \theta_1-\theta^{(1)}_1\\\theta_2-\theta^{(1,2)}_2
  \end{array}
\right)\,.
\]
By \eqref{eq:52}, the gradient is zero at 
$(\theta_1^{(1,2,3)},\theta_2^{(1,2,3)})$ so that, due to the function
$\theta_3(\theta_1,\theta_2)$ being strictly concave,
$(\theta_1^{(1,2,3)},\theta_2^{(1,2,3)})=(\theta_1^{(1)},\theta_2^{(1,2)})$\,.
\end{proof}
\begin{remark}
  By \eqref{eq:37} and  conservation of momenta, the cost to get to  $x$ is
\[  V^{(1),(1,2)}(x)=\theta^{(1)}_1x_1+\theta^{(1,2)}_2x_2+\theta^{(1,2,3)}_3
x_3=(\theta^{(1)}_1,\theta^{(1,2)}_2,\theta^{(1,2,3)}_3)x=\theta^{(1,2,3)}x\,.
\]
where the superscript at $V$ indicates that this cost is associated
with travelling down the $x_1$--axis first and in the $x_1x_2$--plane next.
In order to find the actual optimal trajectory, we could proceed in reverse
time: first, find $s^{(1,2)}=(s^{(1,2)}_1,s^{(1,2)}_2,0)$ 
such that $x-s^{(1,2)}$ is colinear
to $\nabla H_0(\theta^{(1,2,3)})$\,, then find
$s^{(1)}=(s^{(1)}_1,0,0)$ such that the slope of
$s^{(1,2)}-s^{(1)}$ belongs to $\partial H_{3}(\theta^{(1,2)})$\,.
The polygonal curve with corner points at $O$\,, $s^{(1)}$\,, $s^{(1,2)}$
and $x$ is the optimal trajectory.
\end{remark}
\begin{remark}
  The proof implies that under the hypotheses 
$\theta_1^{(1)}<\theta_1(x)$ and 
$\theta_2^{(1,2)}<\theta_2(x)$\,.
\end{remark}
\begin{remark}
   If $s^{(1)}_1=0$\,, then 
  an optimal trajectory may be assumed to follow a straight line 
 in the interior of the
 $x_1x_2$--quarter plane.
 Then, $p^{(1,2)}=s^{(1,2)}$\,, the piece of the 
 optimal trajectory in the $x_1x_2$--plane is a straight line
 segment that connects the origin to  $s^{(1,2)}$ so that
 $s^{(1,2)}_1>0$ and $s^{(1,2)}_2>0$\,,
   and $\theta^{(1,2)}_1=\theta_1(s^{(1,2)})$\,,
with $\theta_1(s^{(1,2)})$ being defined by \eqref{eq:22'}.
 We still have \eqref{eq:35}, \eqref{eq:36} and the equation 
$   \dot q(t)\in\partial H_3(\theta^{(1,2)})\,.$
Moreover, $\theta^{(1,2)}= \theta(s^{(1,2)})$\,, with $\theta(s^{(1,2)})$
  satisfying the equation
$H_0(\theta(s^{(1,2)}))=0$\,. 
  The cost of going down
   $s^{(1,2)}$ is $\theta(s^{(1,2)}) s^{(1,2)}$\,.
 The cost of getting to  point $x\in\R_+^3$ is, in analogy with
 \eqref{eq:37}, 
$\theta(s^{(1,2)})s^{(1,2)}+
 \theta(x-s^{(1,2)})(x-s^{(1,2)})\,,
  $
where
 $\theta_1(x-s^{(1,2)})=\theta_1(s^{(1,2)})$ and
 $\theta_2(x-s^{(1,2)})=\theta_2(s^{(1,2)})$\,.
 In addition,
 $H_0(x-s^{(1,2)})=0$\,.
 Hence, the cost of getting to $x$ by going in the interior of the
 $x_1x_2$--plane first  and next taking off to $x$ is 
$     V^{(1,2)}(x)=\inf_{s^{(1,2)}}(\theta_1(s^{(1,2)})x_1+
\theta_2(s^{(1,2)})x_2+\theta_3(x-s^{(1,2)})x_3)\,.$
\end{remark}
Putting everything together yields the following prescription for the
optimal cost of getting to $x$\,
\[V(x)=V^{(1),(1,2)}(x)\wedge
V^{(1),(1,3)}(x)\wedge V^{(2),(2,3)}(x)\wedge V^{(2),(2,1)}(x)
\wedge V^{(3),(3,1)}(x)\wedge V^{(3),(3,2)}(x)\wedge
V^{(1,2,3)}(x)\,,\]
where $V^{(1,2,3)}(x)$ stands for the cost of getting to $x$ in a
straight line, i.e. $V^{(1,2,3)}(x)=\theta(x)x$\,. 

\section{Some generalities}
\label{sec:some-generalities}

In this section, we give a general characterisation of optimal
trajectories. It summarises the developments above and  extends
Puhalskii \cite{Puh19}.

\begin{lemma}
  \label{le:optim}
Suppose that  trajectory $q$  gets from the origin to $x$ by traversing
 facets $F_{J_1},\ldots,F_{J_k}$ such that
 $J_1\supset\ldots\supset J_k$
with respective  momenta $\theta^{(J_1)},\ldots, \theta^{(J_k)}$\,, i.e.
$\dot q(t)\in\partial H_{J_l}(\theta^{(J_l)})$ a.e. when 
$q(t)\in F_{J_l}$ and $H_{J_l}(\theta^{(J_l)})=0$\,, where $l\in\{1,\ldots,k\}$\,. 
Suppose that
 the
following conservation law holds: 
$\theta^{(J_{l})}_i=\theta^{(J_{l-1})}_i$ provided  $i\notin J_{l-1}$\,.
Then the  trajectory's cost is $\theta^{(J_{k})}x$\,.
\end{lemma}
\begin{proof}
If $i\notin J_l$\,, then
the conservation law implies that
$\theta^{(J_l)}_i=\theta^{(J_{l+1})}_i=\ldots=\theta^{(J_k)}_i$\,.
Let $t_k$ denote the length of time it takes to get to the
destination. Then
the associated cost is
\begin{multline*}
  \int_0^{t_k}L(q(s),\dot q(s))\,ds=
\sum_{l=1}^k \int_0^{t_k}\ind_{F_{J_l}}(q(s))L_{J_l}(\dot q(s))\,ds=
\sum_{l=1}^k \int_0^{t_k}\ind_{F_{J_l}}(q(s))(\theta^{(J_l)}\dot
q(s)-
H_{J_l}(\theta^{(J_l)}))\,ds\\=
\sum_{l=1}^k \int_0^{t_k}\ind_{F_{J_l}}(q(s))\theta^{(J_l)}\dot
q(s)\,ds=\sum_{l=1}^k \int_0^{t_k}\ind_{F_{J_l}}(q(s))\theta^{(J_k)}\dot
q(s)\,ds=\theta^{(J_k)}q(t_k)=\theta^{(J_k)}x\,.
\end{multline*}
\end{proof}
We conjecture the conservation law in the statement of the lemma
to hold generally.
Under the conservation law, in order to find $\theta^{(J_{l})}_i$\,, for $i\in J_{l-1}$\,, 
given $\theta^{(J_{l-1})}$\,, the
following equations need to be solved:
$h_{j}(\theta^{(J_{l})})=0\,, j\in J_{l}\,,$ and
$H_0(\theta^{(J_{l})})=0\,. $
We thus have $\abs{J_{l-1}}$ unknowns and $\abs{J_l}+1$ equations.
(Because of the conservation law,
  $K-\abs{J_{l-1}}$ components carry over from $\theta^{(J_{l-1})}$\,.)
If the number of zero components decreases by one when transiting from 
$F_{J_{l-1}}$ to $F_{J_l}$\,, then the unknowns could be generally
found.
\appendix

\section{Proof of Theorem \ref{the:hamilt}}
\label{sec:appendix}
\begin{lemma}
  \label{le:alpha}
  \begin{enumerate}
  \item 
The functions $f_k(u)$ and $g_k(u)$ are well defined, are infinitely 
 differentiable, with  positive first
derivatives,
 and  range in  $(-\infty,\infty)$ 
 as $u$ ranges in $(-1,\infty)$\,.
In addition,  $f_k(0)=g_k(0)=0$ and 
$f_k'(0)=g_k'(0)=1$\,.
\item
The functions $f_k(e^\vartheta-1)-\vartheta$  are
nondecreasing for $\vartheta>0$\,. 
The functions $f_k(e^\vartheta-1)$  are convex,
being strictly convex provided the $\xi_k$ 
  are nonconstant.
\item The functions $g_k(h_k(\theta))$ and $g_k(h_k(\theta)^+)$ are convex.
  \end{enumerate} \end{lemma}

 \begin{proof}

   The function $f_k$ is well defined, continuous, increasing
 and  ranges from $-\infty$ to
$\infty$ as $\vartheta$ ranges from $-\infty$ to $\infty$ 
because $\mathbf
P(\xi_k=0)=0$\,. Furthermore, the dominated convergence theorem and
the implicit function theorem imply that the function is
infinitely  differentiable. 
It is readily seen that  $f_k(0)=0$\,.
Let $\tilde f_k(\vartheta)=f_k(e^{\vartheta}-1)$\,.
 By Jensen's inequality,
 $\tilde f_k(\vartheta)\ge \vartheta$\,.
Exponentiating and differentiating 
in \eqref{eq:1a} yield
\begin{equation}
  \label{eq:7a}
 \tilde f_k'(\vartheta)\mathbf E e^{-
 \tilde f_k(\vartheta)\lambda_k\xi_k}\lambda_k\xi_k
=e^{-\vartheta}\,.
\end{equation}
Hence, $
 \tilde f_k'(0)=1$ and $\tilde f'_k(\vartheta)>0$\,. 
We differentiate \eqref{eq:7a}:
\begin{equation}
  \label{eq:8}
\tilde   f_k''(\vartheta)\mathbf E e^{-
  \tilde f_k(\vartheta)\lambda_k\xi_k}\lambda_k\xi_k=
( \tilde f_k'(\vartheta))^2\mathbf E e^{-
  \tilde f_k(\vartheta)\lambda_k\xi_k}(\lambda_k\xi_k)^2
-e^{-\vartheta}\,.
\end{equation}
By the Cauchy--Schwarz inequality and \eqref{eq:7a}, 
\[
   \tilde f_k'(\vartheta)^2\mathbf E e^{-
  \tilde f_k(\vartheta)\lambda_k\xi_k}(\lambda_k\xi_k)^2
\ge \bl( \tilde f_k'(\vartheta)\mathbf E e^{-
  \tilde f_k(\vartheta)\lambda_k\xi_k}(\lambda_k\xi_k)\br)^2
/\mathbf E e^{-
  \tilde f_k(\vartheta)\lambda_k\xi_k}=e^{-\vartheta}\,,
\]
the inequality being strict when $\xi_k$ is nonconstant.
By \eqref{eq:8}, $ \tilde f_k''(\vartheta)\ge0$, so, $
\tilde f_k'(\vartheta)$ is nondecreasing
and $\tilde f_k(\vartheta)$ is convex. Hence, $ \tilde f_k'(\vartheta)\ge1$ when $\vartheta>0$\,,
so,
$\tilde f_k(\vartheta)-\vartheta$ is nondecreasing in $\vartheta>0$\,.

We prove part 3.
Let, given $\hat K\le K$\,, 
\begin{equation}
  \label{eq:25}
  \hat h_k(\theta)=\ln\bl(\sum_{l=1}^{\hat K}(e^{\theta_l}-1)p_{kl}+1\br)\,.
\end{equation}
We show it is a convex function. The proof is done iteratively: for
$\hat K=1\,, \hat K=2\,, \ldots,\hat K=K$\,. For
$\hat K=1$\,, a direct computation shows that the second derivarive of $\hat
h_k$ is nonnegative. Let $\hat K=2$ so that
$\theta=(\theta_1,\theta_2)$ and 
\begin{equation}
  \label{eq:28}
    \hat h_k(\theta)=\ln\bl((e^{\theta_1}-1)p_{k1}
+(e^{\theta_2}-1)p_{k2}+1\br)\,.
\end{equation}
We have that
\begin{equation}
  \label{eq:50}
   \nabla^2\hat h_k(\theta)=
D^{-2}\bl(e^{\theta_1}p_{k1}e^{\theta_2}p_{k2}
\left(\begin{array}{cc}
  1&
-1\\
- 1&1
\end{array}\right)
+(1-p_{k1}-p_{k2})
\left(\begin{array}{cc}
  e^{\theta_1}p_{k1}&
0\\
0&e^{\theta_2}p_{k2}
\end{array}\right)\br)\,,
\end{equation}
where
\[  D=(e^{\theta_1}-1)p_{k1}
+(e^{\theta_2}-1)p_{k2}+1\,.
\]
All principal minors of the matrices on the  righthand side of \eqref{eq:50}
 being nonnegative, the matrix $\nabla^2\hat
h_k(\theta)$ is positive semidefinite, so, $\hat h_k$ is convex.
 
 Suppose that $\hat h_k$ is convex for some $\hat K\le K-1$
and let $\tilde h_k$ represent $\hat h_k$  orresponding to $\hat
K+1$\,, i.e.
\begin{equation}
  \label{eq:26}
  \tilde h_k(\theta)=\ln\bl(\sum_{l=1}^{\hat K+1}(e^{\theta_l}-1)p_{kl}+1\br)\,.
\end{equation}
By \eqref{eq:25} and \eqref{eq:26},
\begin{equation}
  \label{eq:27}
  \tilde h_k(\theta)=\ln(e^{\hat h_k(\theta)}+(e^{\theta_{\hat
      K+1}}-1)p_{k,\hat K+1})\,.
\end{equation}
Since the righthand side of \eqref{eq:28} is a convex function of
$(\theta_1,\theta_2)$\,,   $\tilde h_k(\theta)$ is a convex function
of $\hat h_k(\theta)$ and $\theta_{\hat K+1}$\,. By  hypothesis, 
$\hat h_k(\theta)$ is a convex function of
$\theta_1,\ldots,\theta_{\hat K}$\,. Hence, $\tilde h_k(\theta)$ is a convex
function of $\theta_1,\ldots,\theta_{\hat K+1}$\,, as needed.

By \eqref{eq:4'} and \eqref{eq:25} with $\hat K=K$\,,
$g_k(h_k(\theta))=g_k(e^{-\theta_k}e^{\hat h_k(\theta)}-1)$\,. 
By  part
1, 
$g_k(h_k(\theta))$ is a convex function of $-\theta_k+
\hat h_k(\theta)$\,.  Since $\hat h_k(\theta)$ is a convex function of
$\theta$\,,
$g_k(h_k(\theta))$ is a convex function of $\theta$ too.
A similar proof works for $g_k(h_k(\theta)^+)$\,.
 \end{proof}

\begin{proof}[Proof of Theorem \ref{the:hamilt}]
  
The proof draws on 
 Puhalskii \cite{Puh19}.
By the definition of $f_k$\,, for $\vartheta\in\R$ and $t\in\R$\,,
$\ln 
\mathbf Ee^{\vartheta \xi_k}=-t$ if and only if
$\vartheta=-\lambda_kf_k(e^t-1)$\,. Hence,
$\sup_{\vartheta \in\R}(\vartheta-a_k\ln\mathbf Ee^{\vartheta \xi_k})=
\sup_{\vartheta\in\R}(a_k\vartheta-\lambda_kf_k(e^\vartheta-1))$\,.
Similarly, $\sup_{\vartheta \in\R}(\vartheta-d_k\ln\mathbf Ee^{\vartheta \eta_k})=
\sup_{\vartheta\in\R}(d_k\vartheta-\mu_kg_k(e^\vartheta-1))$\,.
Also, $  \pi(x)=\sup_{\vartheta\in\R}(\vartheta x-(e^\vartheta-1))\,,x\in\R\,.$
Owing to  \eqref{eq:psi_J}, \eqref{eq:Psi}, and the minimax theorem,
assuming that $y=(y_k)$\,, $a=(a_k)$\,, $d=(d_k)$\,, and $r=(r_{kl})$\,,
\begin{multline*}
L_J(y)=  \inf_{\substack{(a,d,r)\in\R_+^K
\times\R_+^K\times\mathbb{S}_+^{K\times K}
:\\y=a+(r^T-I)d}}\psi_J(a,d,r)
\\= \inf_{\substack{(a,d,r)\in\R_+^K
\times\R_+^K\times\mathbb{S}_+^{K\times K}
:\\y=a+(r^T-I)d}}\Bl(\sum_{k=1}^K
\sup_{\vartheta\in\R}\bl(\vartheta\,a_k-f_k(e^\vartheta-1)\lambda_k\br)\,
+\sum_{k\in J^c}
\sup_{\vartheta\in\R}\bl(\vartheta\,d_k-g_k(e^\vartheta-1)\mu_k\br)\,
\\
+\sum_{k\in
  J}\sup_{\vartheta\in\R}\bl(\vartheta\,d_k-g_k(e^\vartheta-1)\mu_k\br)\,
\,
\ind_{\{d_k>\mu_k\}}
+\sum_{k=1}^K d_k \biggl[\sum_{l=1}^K 
\sup_{\vartheta\in\R}\bl(\vartheta\,r_{kl}-(e^\vartheta-1)p_{kl}\br)\,
\\+
\sup_{\vartheta\in\R}\bl(\vartheta\,\bl(1-\sum_{l=1}^K
  r_{kl}\br)-(e^\vartheta-1)
\bl(1-\sum_{l=1}^K
p_{kl}\br)\br)\,\biggr]\Br)\\=
 \inf_{\substack{(a,d,r)\in\R_+^K
\times\R_+^K\times\mathbb{S}_+^{K\times K}
:\\y=a+(r^T-I)d}}\,\sup_{\substack{(\theta_k),(\vartheta_k),\\(\sigma_{kl}),(\rho_k)}}
\Bl(\sum_{k=1}^K
\bl(\theta_k\,a_k-f_k(e^{\theta_k}-1)\lambda_k\br)\,
+\sum_{k\in J^c}
\bl(\vartheta_k\,d_k-g_k(e^{\vartheta_k}-1)\mu_k\br)\,
\\
+\sum_{k\in
  J}\bl(\vartheta_k\,d_k-g_k(e^{\vartheta_k}-1)\mu_k\br)\,
\,
\ind_{\{d_k>\mu_k\}}
+\sum_{k=1}^K d_k \biggl[\sum_{l=1}^K 
\bl(\sigma_{kl}\,r_{kl}-(e^{\sigma_{kl}}-1)p_{kl}\br)\,
\\+
\rho_k\,\bl(1-\sum_{l=1}^K
  r_{kl}\br)-(e^{\rho_k}-1)
\bl(1-\sum_{l=1}^K
p_{kl}\br)\,\biggr]\Br)\\=
\,\sup_{\substack{(\theta_k),(\vartheta_k),\\(\sigma_{kl}),(\rho_k)}}
 \inf_{\substack{(a,d,r)\in\R_+^K
\times\R_+^K\times\mathbb{S}_+^{K\times K}
:\\y=a+(r^T-I)d}}
\Bl(\sum_{k=1}^K
\bl(\theta_k\,a_k-f_k(e^{\theta_k}-1)\lambda_k\br)\,
+\sum_{k\in J^c}
\bl(\vartheta_k\,d_k-g_k(e^{\vartheta_k}-1)\mu_k\br)\,
\\
+\sum_{k\in
  J}\bl(\vartheta_k\,d_k-g_k(e^{\vartheta_k}-1)\mu_k\br)\,
\,
\ind_{\{d_k>\mu_k\}}
+\sum_{k=1}^K d_k \biggl[\sum_{l=1}^K 
\bl(\sigma_{kl}\,r_{kl}-(e^{\sigma_{kl}}-1)p_{kl}\br)\,
\\+
\rho_k\,\bl(1-\sum_{l=1}^K
  r_{kl}\br)-(e^{\rho_k}-1)
\bl(1-\sum_{l=1}^K
p_{kl}\br)\,\biggr]\Br)\\=
\,\sup_{\substack{(\theta_k),(\vartheta_k),\\(\sigma_{kl}),(\rho_k)}}
 \inf_{\substack{(a,d,r)\in\R_+^K
\times\R_+^K\times\mathbb{S}_+^{K\times K}
:\\y=a+(r^T-I)d}}
\Bl(\sum_{k=1}^K
\theta_k\,a_k
+\sum_{k\in J^c}
\vartheta_k\,d_k
\\
+\sum_{k\in
  J}\vartheta_k\,d_k\,
\ind_{\{d_k>\mu_k\}}
+\sum_{k=1}^K d_k \sum_{l=1}^K 
\sigma_{kl}\,r_{kl}
+\sum_{k=1}^K d_k
\rho_k\,\bl(1-\sum_{l=1}^K
  r_{kl}\br)\\
-\sum_{k=1}^Kf_k(e^{\theta_k}-1)\lambda_k
-\sum_{k\in J^c}g_k(e^{\vartheta_k}-1)\mu_k\,
-\sum_{k\in
  J}g_k(e^{\vartheta_k}-1)\mu_k
\,
\ind_{\{d_k>\mu_k\}}\\
-\sum_{k=1}^K d_k \sum_{l=1}^K (e^{\sigma_{kl}}-1)p_{kl}
-\sum_{k=1}^K d_k(e^{\rho_k}-1)
\bl(1-\sum_{l=1}^K
p_{kl}\br)
\Br)\\=
\,\sup_{\substack{(\theta_k),(\vartheta_k),\\(\sigma_{kl}),(\rho_k)}}
 \inf_{\substack{(d,r)\in\R_+^K\times\mathbb{S}_+^{K\times K}
:\\y\ge(r^T-I)d}}
\Bl(\sum_{k=1}^K
\theta_k\,(y_k+d_k-\sum_{l=1}^Kr_{lk}d_l)
+\sum_{k\in J^c}
\vartheta_k\,d_k
\\
+\sum_{k\in
  J}\vartheta_k\,d_k\,
\ind_{\{d_k>\mu_k\}}
+\sum_{k=1}^K d_k \sum_{l=1}^K 
\sigma_{kl}\,r_{kl}
+\sum_{k=1}^K d_k
\rho_k\,\bl(1-\sum_{l=1}^K
  r_{kl}\br)\\
-\sum_{k=1}^Kf_k(e^{\theta_k}-1)\lambda_k
-\sum_{k\in J^c}g_k(e^{\vartheta_k}-1)\mu_k\,
-\sum_{k\in
  J}g_k(e^{\vartheta_k}-1)\mu_k
\,
\ind_{\{d_k>\mu_k\}}\\
-\sum_{k=1}^K d_k \sum_{l=1}^K (e^{\sigma_{kl}}-1)p_{kl}
-\sum_{k=1}^K d_k(e^{\rho_k}-1)
\bl(1-\sum_{l=1}^K
p_{kl}\br)
\Br)\,.
\end{multline*}
Let us minimise with respect to $r=(r_{kl})$ the sum
\[
  \sum_{k=1 }^K\bl(-\sum_{l=1}^K
\theta_l r_{kl}+\sum_{l=1}^K 
\sigma_{kl}\,r_{kl}-\rho_k\sum_{l=1}^K
  r_{kl}
\br)\,d_k\,
\]
subject to the constraint that $r_{kl}\ge0$\,,
 $\sum_k r_{kl}d_k\le y_l+d_l$ and $\sum_lr_{kl}\le1$\,.
By  Lagrange multipliers,
there exists $K$-tuple 
$\alpha=(\alpha_1,\ldots,\alpha_K)$ with nonnegative entries
such that 
\[
\sum_{k=1 }^K\bl(-\sum_{l=1}^K
\theta_l r_{kl}+\sum_{l=1}^K 
\sigma_{kl}\,r_{kl}-\rho_k\sum_{l=1}^K
  r_{kl}
\br)\,d_k\,+\sum_{l=1}^K\alpha_l \sum_{k=1}^K r_{kl}d_k
\]
attains  minimum given that
 $r_{kl}\ge0$ and $\sum_lr_{kl}\le1$\,. 
The expression to be minimised can be written as
\[  \sum_{k=1 }^K\sum_{l=1}^K\bl(-
\theta_l +
\sigma_{kl}-\rho_k
+ \alpha_l 
\br)\,d_kr_{kl}\,.
\]
The minimum   equals
$  \sum_{k }d_k\min_{1\le l \le K}\bl(-
\theta_l +
\sigma_{kl}-\rho_k
+ \alpha_l 
\br)\wedge 0$\,.
It is attained at $r$ defined as follows. If $
\min_{1\le l \le K}\bl(-
\theta_l +
\sigma_{kl}-\rho_k
+ \alpha_l 
\br)<0$ and  the latter minimum is
attained at $m$\, then $r_{km}=1$ and $r_{kl}=0\,, l\not=m$.
 If $\min_{1\le l \le K}\bl(-
\theta_l +
\sigma_{kl}-\rho_k
+ \alpha_l 
\br)\ge0$\,, then 
 $r_{kl}=0$\,, for all $l$\,.

The maximum over $\alpha_l$ is attained for  $\alpha_l=0$\,. 
It follows that
\begin{multline}
  \label{eq:16}
  L_J(y)
 = \,\sup_{\substack{(\theta_k),(\vartheta_k),\\(\sigma_{kl}),(\rho_k)}}
 \inf_{(d_k)}
\Bl(\sum_{k=1}^K
\theta_k\,y_k+
  \sum_{k=1 }^K(c_k+\theta_k+\vartheta_k)d_k-\sum_{k\in J}\vartheta_k\ind_{\{d_k\le \mu_k\}}d_k
-\sum_{k=1}^Kf_k(e^{\theta_k}-1)\lambda_k
\\-\sum_{k\in J^c}g_k(e^{\vartheta_k}-1)\mu_k\,
-\sum_{k\in
  J}g_k(e^{\vartheta_k}-1)\mu_k
\,
\ind_{\{d_k>\mu_k\}}\Br)\,,
\end{multline}
where
\begin{equation}
  \label{eq:19}
  c_k=
 \min_{1\le l \le K}(-\theta_l +
\sigma_{kl}-\rho_k)\wedge0+
\rho_k
-\sum_{l=1}^K (e^{\sigma_{kl}}-1)p_{kl}
- (e^{\rho_k}-1)
\bl(1-\sum_{l=1}^K
p_{kl}\br)\,.
\end{equation}

If $c_k+\theta_k+\vartheta_k<0$\,, for some $k$\,, 
then letting $d_k\to\infty$ shows that the $\inf $
in \eqref{eq:16} equals $-\infty$\,.  Therefore,
\begin{multline*}
L_J(y)=  \,\sup_{\substack{(\theta_k),\,(\vartheta_k),\,
(\sigma_{kl}),\,(\rho_k):\\
c_k+\theta_k+\vartheta_k\ge0}} \inf_{d_k\ge0}
\Bl(\sum_{k=1}^K
\theta_k\,y_k+
  \sum_{k=1 }^K(c_k+\theta_k+\vartheta_k)d_k-\sum_{k\in J}\vartheta_k\ind_{\{d_k\le \mu_k\}}d_k\\
-\sum_{k=1}^Kf_k(e^{\theta_k}-1)\lambda_k
-\sum_{k\in J^c}g_k(e^{\vartheta_k}-1)\mu_k\,
-\sum_{k\in
  J}g_k(e^{\vartheta_k}-1)\mu_k
\,
\ind_{\{d_k>\mu_k\}}\Br)\,.
\end{multline*}
In the sum, it is optimal to take $d_k=0$
for $k\in J^c$ so that
\begin{multline*}
L_J(y)=  \sup_{\substack{(\theta_k),\,(\vartheta_k),\,
(\sigma_{kl}),\,(\rho_k):\\
c_k+\theta_k+\vartheta_k\ge0}} \inf_{d_k\ge0}
\Bl(\sum_{k=1}^K
\theta_k\,y_k+  \sum_{k\in J }\bl((c_k+\theta_k+\vartheta_k\ind_{\{d_k> \mu_k\}})d_k
-g_k(e^{\vartheta_k}-1)\mu_k\ind_{\{d_k>\mu_k\}}\br)\\
-\sum_{k=1}^Kf_k(e^{\theta_k}-1)\lambda_k
-\sum_{k\in J^c}g_k(e^{\vartheta_k}-1)\mu_k\Br)\,.
\end{multline*}
By case analysis, recalling that $g_k(e^{\vartheta_k}-1)\ge \vartheta_k$\,,
\[
\inf_{d_k\ge0}\bl((c_k+\theta_k+\vartheta_k\ind_{\{d_k> \mu_k\}})d_k
-g_k(e^{\vartheta_k}-1)\mu_k\ind_{\{d_k>\mu_k\}}\br)=
\mu_k(c_k+\theta_k+\vartheta_k-g_k(e^{\vartheta_k}-1))\wedge0\,.
\]
Consequently, 
\begin{multline*}
  L_J(y)=\sup_{\substack{(\theta_k),\,(\vartheta_k),\,
(\sigma_{kl}),\,(\rho_k):\\
c_k+\theta_k+\vartheta_k\ge0}} \Bl(\sum_{k=1}^K
\theta_k\,y_k+\sum_{k\in J }\mu_k(c_k+\theta_k+\vartheta_k-
g_k(e^{\vartheta_k}-1))\wedge0\\-\sum_{k=1}^Kf_k(e^{\theta_k}-1)\lambda_k
-\sum_{k\in J^c}g_k(e^{\vartheta_k}-1)\mu_k\Br)\,.
\end{multline*}
We now maximise the $c_k$ in \eqref{eq:19} over  $\sigma_{kl}$ and $\rho_k$\,.
It seems  easier to minimise 
\begin{equation}
  \label{eq:20}
 \sum_{l=1}^K (e^{-\sigma_{kl}}-1)p_{kl}
+(e^{-\rho_k}-1)
\bl(1-\sum_{l=1}^K
p_{kl}\br)+
 \max_{1\le l\le K}(\theta_l +
\sigma_{kl})\vee\rho_k\,,
\end{equation}
which is
 a convex function of $(\sigma_{kl}\vert_{l=1}^K,\rho_k)$\,.  
Let 
 $\beta_k=
\max_{1\le l\le K}(\theta_l +
\sigma_{kl})\vee\rho_k$\,.
 The subdifferential of the function in \eqref{eq:20} is 
\[
 \sum_{l=1}^K (-e^{-\sigma_{kl}}p_{kl})(1_l,0)-
e^{-\rho_k}(1-\sum_{l=1}^Kp_{kl})(0_K,1)+
\text{co}\,
(\sum_{l\in I_K}(1_l,0)
+(0_K,1)\ind_{\{\rho_k=\beta_k\}})\,,
\]
 where $1_l$ represents the $K$-vector, whose $l$-th entry equals $1$\,,
 the rest of the entries being equal to $0$\,,
$0_K$ represents the zero
 $K$--vector,   co stands for taking 
the convex hull of a set of vectors, and
 $I_k$ represents
the set (possibly, null) of
 $l$ such that
$\theta_l+\sigma_{kl}=\beta_k$\,. As the subdifferential at a 
point where the minimum is attained
 must contain the zero vector, it is necessary that
 $I_k$ contain all $l$ such that  $p_{kl}\not=0$ and that
$\rho_k=\beta_k$ provided $\sum_lp_{kl}<1$\,. According to the definition of a
convex hull, there exist nonnegative weights 
 $\alpha_1,\,\ldots,\alpha_K,\alpha_{K+1}$
adding up to 1 such that
$-e^{-\sigma_{kl}}p_{kl}+\alpha_l=0$ for  $l\in I_k$\,, 
$\alpha_l=0$ for $l\not\in I_k$\,,
$-e^{-\rho_k}(1-\sum_lp_{kl})+\alpha_{K+1}=0$ provided
$\rho_k=\beta_k$ and  $\alpha_{K+1}=0$\,, otherwise. Thus,
  $\alpha_l=e^{-\sigma_{kl}}p_{kl}$\,, for all $l$\,, and
$\alpha_{K+1}=e^{-\rho_k}(1-\sum_lp_{kl})$\,.
Since the $\alpha_i$ sum to 1,
$\sum_{l}e^{-\sigma_{kl}}p_{kl}+e^{-\rho_k}(1-\sum_lp_{kl})=1$\,. 
Since $\theta_l+\sigma_{kl}=\beta_k$ provided $p_{kl}\not=0$ and
 $\rho_k=\beta_k$ provided $1-\sum_lp_{kl}>0$\,,
$\sum_{l}e^{-\beta_k+\theta_l}p_{kl}+e^{-\beta_k}(1-\sum_lp_{kl})=1$\,,
i.e.
the optimal value of $\beta_k$ is given by the relation
\begin{equation}
e^{\beta_k}=\sum_{l=1}^K(e^{\theta_l}-1)p_{kl}+1\,. \label{eq:7'}
\end{equation}
The minimum of the function in \eqref{eq:20} equals
$\sum_l(e^{\theta_l-\beta_k}-1)p_{kl}+(e^{-\beta_k}-1)(1-\sum_lp_{kl})+\beta_k=
e^{-\beta_k}(\sum_l(e^{\theta_l}-1)p_{kl}
+1)-1+\beta_k=
\beta_k$\,.

Thus, 
  \begin{multline}
\label{eq:6a}
L_J(y) =\sup_{(\theta_k),\,(\vartheta_k):\,
\theta_k+
\vartheta_k\ge\beta_k}
\Bl(\sum_{k=1}^K
\theta_k\,y_k
+\sum_{k\in J}\mu_k\bl(-\beta_k+\theta_k
+\vartheta_k-g_k(e^{\vartheta_k}-1)\br)\wedge0\\
  -\sum_{k=1}^Kf_k(e^{\theta_k}-1)\lambda_k-\sum_{k\in J^c}g_k(e^{\vartheta_k}-1)\mu_k\,
\Br)\,.
\end{multline}
 If $-\beta_k+\theta_k\ge0$\,, for some $k\in J$\,,
then recalling that $\vartheta_k-g_k(e^{\vartheta_k}-1)\le0$ and
$g_k(0)=0$\,,
we can take  $\vartheta_k=0$\,. It follows that
the righthand side of \eqref{eq:6a}
 can be equivalently written as
\begin{multline*}
\,\sup_{(\theta_k),\,(\vartheta_k):\,
\theta_k+
\vartheta_k\ge\beta_k}
\Bl(\sum_{k=1}^K
\theta_k\,y_k
+\sum_{k\in J}\ind_{\{\theta_k<\beta_k\}}\mu_k\bl(-\beta_k+\theta_k
+\vartheta_k-g_k(e^{\vartheta_k}-1)\br)\wedge0\,\\
-\sum_{k=1}^Kf_k(e^{\theta_k}-1)\lambda_k-\sum_{k\in J^c}g_k(e^{\vartheta_k}-1)\mu_k\,
\Br)\,.
\end{multline*}
As the function $\vartheta-g_k(e^\vartheta-1)$ 
is nondecreasing for  $\vartheta\ge0$\,, the supremum over
$\vartheta_k$ for $k\in J$ is attained on the boundary of the set of
constraints, i.e. $\vartheta_k=\beta_k-\theta_k$. We obtain that
\begin{multline*}
L_J(y)=\sup_{(\theta_k)}
\Bl(\sum_{k=1}^K
\theta_k\,y_k
-\sum_{k\in J}g_k((e^{\beta_k-\theta_k}-1)^+)\mu_k
-\sum_{k=1}^Kf_k(e^{\theta_k}-1)\lambda_k-\sum_{k\in
  J^c}g_k(e^{\beta_k-
\theta_k}-1)\mu_k\,
\Br)\,.
\end{multline*}
 Recalling 
\eqref{eq:4'} and \eqref{eq:7'}
implies that $e^{\beta_k-\theta_k}
=1+h_k(\theta)$ so that
\[
L_J(y)=\sup_{(\theta_k)}
\Bl(\sum_{k=1}^K
\theta_k\,y_k
-\sum_{k\in J}g_k(h_k(\theta)^+)\mu_k
-\sum_{k=1}^Kf_k(e^{\theta_k}-1)\lambda_k-\sum_{k\in
  J^c}g_k(h_k(\theta))\mu_k\,
\Br)\,.
\]
Recalling \eqref{eq:6} and Lemma \ref{le:alpha} yields \eqref{eq:9}.
\end{proof}
\begin{corollary}
  The function $L_\emptyset(y)$ is strictly convex.
  \begin{proof}
    The function $H_0(\theta)$ is differentiable, so $L_\emptyset(y)$ is
    strictly convex as a convex conjugate of $H_0(\theta)$\,.
  \end{proof}
\end{corollary}
\begin{figure}[h!]
  \centering
  \includegraphics[scale=.7]{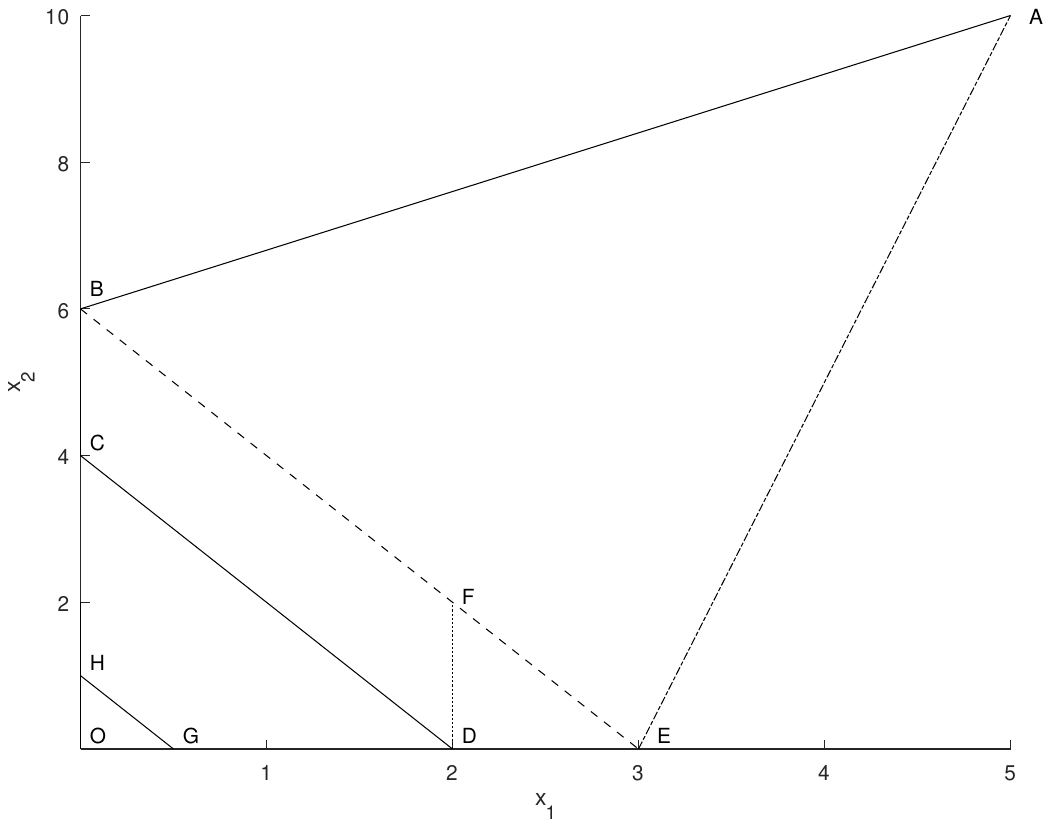}
  \caption{Optimality of a two--straight--line--segment path}
      \label{fig:6}
\end{figure}

\section{Optimality of  two straight--line segments in 2D}
\label{sec:optim-two-stra}
Let $\hat q$
represent an optimal trajectory to get from the origin to $x$ 
and let $\hat t$ represent the length of  time it takes to reach the destination. 

We have that $\hat q(\hat t)=x$\,. 
It is useful to
introduce the time reversed process $\overline q(t)=\hat q(\hat
t-t)$\,. 
The trajectory $\overline q$ consists of
 two straight--line segments  at most. This follows from geometric
 considerations, as in Fig. \ref{fig:6}.
If one assumes that the path $ABCDO$ is optimal,  then the cost of $CDO$ 
is less than  or equal to the cost of 
  $CO$\,. By similarity, the cost of $FED$ is less than or is equal to
  the cost
$BC$\,. Since $BF$ and $CD$ have the same cost, by translation invariance,
the  path $ABEO$ is also optimal, which is impossible as, by
srict convexity, the cost of path $AE$ is less than the cost of 
path $ABE$\,. Hence, $ABCDO$ cannot be optimal. If the path $ABCDO$
leaves the horizontal axis  before hitting the origin at
some point $G$ and hits the vertical axis at some point $H$\,,
then the path $ABCDGH$ cannot be optimal for a similar reason.
As a consequence, an optimal path to get from the origin to a given
point is either the straight--line segment connecting the origin and
the destination or a path that 
consists  of a straight--line segment along one of the axes
and a straight--line segment to the destination. 

\section{Motion down an axis}
\label{sec:optimal-paths}

\begin{lemma}
  Suppose that the network is ergodic. Given $m\in \mathcal{K}$\,,
there exists unique $\theta^{(m)}\not=0$ such that  
 $h_k(\theta^{(m)})=0$ when 
$k\not=m$ and $H_{\mathcal{K}\setminus\{m\}}(\theta^{(m)})=H_0(\theta^{(m)})=0$\,. Furthermore, 
 $\theta^{(m)}\ge 0$ entrywise, $\theta^{(m)}_m>0$\,,
$e^{\theta^{(m)}_m}\le \rho_m^{-1}$\,, and 
$\partial_mH_0(\theta^{(m)})>0$\,.
\end{lemma}
\begin{proof}
Let $\theta$ satisfy the equations $h_k(\theta)=0$ when $k\not=m$\,.
By  Puhalskii \cite{Puh19},
\begin{equation}
  \label{eq:22}
e^{\theta_l}-1=a_{ml}(e^{\theta_m}-1)\,,  
\end{equation}
 where
 \begin{equation}
   \label{eq:23}
a_{ml}=\frac{c_{ml}}{c_{mm}}\,, 
 \end{equation}
$c_{ml}$ denoting the $(m,l)$--entry of  $C=(I-P^T)^{-1}$\,.
By \eqref{eq:9},
\[
       H_{\mathcal{K}\setminus\{m\}}(\theta)=H_0(\theta)=\sum_{k}f_k(a_{mk}(e^{\theta_m}-1))\lambda_k
+g_m\bl(
h_m(\theta)\br)\mu_m\,.
\]
By \eqref{eq:4'} and \eqref{eq:22},
\[
  h_m(\theta)=(1-e^{-\theta_m})\bl(\sum_{l=1}^Ka_{ml}p_{ml}-1\br)\,.
\]
One can see from \eqref{eq:23}, see also Puhalskii \cite{Puh19}, that
\[
  \sum_{l=1}^Ka_{ml}p_{ml}-1=-\frac{1}{c_{mm}}\,.
\]
We obtain that
\begin{equation}
  \label{eq:24}
         H_0(\theta)=\sum_{k=1}^Kf_k(a_{mk}(e^{\theta_m}-1))\lambda_k
+g_m\bl(c_{mm}^{-1}(e^{-\theta_m}-1)\br)\mu_m
\,.
\end{equation}
The righthand side is the sum of an increasing and a decreasing
function of $\theta_{m}$\,. It is zero at $\theta_{m}=0$\,.
Its derivative at $\theta_m=0$ is
$\sum_{k=1}^Ka_{mk}\lambda_k-c_{mm}^{-1}\mu_m=
(\nu_m-\mu_m)/c_{mm}<0$\,, where $\nu_m$ is the $m$-th entry of
 $\nu=(I-P^T)^{-1}\lambda$ (see Puhalskii \cite{Puh19})\,.
On the other hand, as $\theta_m\to\infty$\,,
$f_m(e^{\theta_m}-1)\to\infty$ and
$g_m(c_{mm}^{-1}(e^{-\theta_m}-1))\to g_m(-c_{mm}^{-1})$\,.
Furthermore, the functions on the righthand side of \eqref{eq:24} are
convex by Lemma \ref{le:alpha} and the facts that $a_{mk}\in[0,1]$ and
$c_{mm}\ge1$\,,
see Puhalskii \cite{Puh19}, so, the derivative of the righthand side is an increasing
function of $\theta_m$\,. It follows that there exists unique
$\theta^{(m)}_m>0$ such that $H_0(\theta^{(m)})=0$\,. In addition, 
$\partial_mH_0(\theta^{(m)})>0$\,. By \eqref{eq:22},
$\theta^{(m)}_k\ge0$\,, for all $k\not=m$\,. 

Since $f_k(\alpha(e^\vartheta-1))\ge\alpha \vartheta$\,, 
$g_k(\alpha(e^\vartheta-1))\ge\alpha \vartheta$
and $H_0(\theta^{(m)})=0$\,, we have, by \eqref{eq:24},  that
$  \sum_{k=1}^Ka_{mk}(e^{\theta^{(m)}_m}-1)\lambda_k
+c_{mm}^{-1}(e^{-\theta^{(m)}_m}-1)\mu_m\le 0\,.
$ Since $\sum_{k=1}^Ka_{mk}\lambda_k=\nu_m/c_{mm}$\,,
$(e^{\theta^{(m)}_m}-1)\nu_m
+ (e^{-\theta^{(m)}_m}-1)\mu_m\le 0$\,, so, $e^{-\theta^{(m)}_m}\ge
\nu_m/\mu_m=\rho_m$\,.  \end{proof}
Let us look at the motion from the origin down the $x_m$--axis, i.e., in the
facet $\Gamma_{\mathcal{K}\setminus\{m\}}=\{x:\,x_m>0\,, x_k=0 \text{ when }k\not=m\}$\,. 
One needs that $\dot q_k(t)=0$\,, for $k\not=m$\,.  
By \eqref{eq:9}, assuming  $\theta=
\theta^{(m)}$\,, we have that
\begin{multline*}
  \dot q_k(t)=\partial_k H_{\mathcal{K}\setminus \{m\}}(\theta^{(m)})=
f'_k(e^{\theta^{(m)}_k}-1)e^{\theta^{(m)}_k}\lambda_k+g'_m\bl(
h_m(\theta^{(m)})\br)\partial_kh_m(\theta^{(m)})\mu_m
\\+\sum_{l\not=m}
\alpha_l\partial_kh_l(\theta^{(m)})\mu_l=
f'_k(e^{\theta^{(m)}_k}-1)e^{\theta^{(m)}_k}\lambda_k+g'_m\bl(
h_m(\theta^{(m)})\br)e^{-\theta^{(m)}_m}e^{\theta^{(m)}_k}p_{mk}\mu_m
\\+\sum_{\substack{l\not=m\,,\\l\not=k}}
\alpha_le^{-\theta^{(m)}_l}e^{\theta^{(m)}_k}p_{lk}\mu_l
-\alpha_ke^{-\theta^{(m)}_k}(1-\sum_{l=1}^Kp_{kl}+\sum_{l\not=k}e^{\theta^{(m)}_l}p_{kl})\mu_k
\,,
\end{multline*}
where $\alpha_l\in[0,1]$\,.
As 
$h_k( \theta^{(m)})=0$\,, 
\[
        1-\sum_{l=1}^Kp_{kl}+\sum_{l\not=k}e^{\theta^{(m)}_l}p_{kl}
=e^{\theta^{(m)}_k}(1-p_{kk})\,.
\]
Hence, we need that, for $k\not=m$\,,
\begin{multline*}
f'_k(e^{\theta^{(m)}_k}-1)\lambda_k+g'_m\bl(
h_m(\theta^{(m)})\br)e^{-\theta^{(m)}_m}p_{mk}\mu_m
+\sum_{l\not=m}
\alpha_le^{-\theta^{(m)}_l}p_{lk}\mu_l
-\alpha_ke^{-\theta^{(m)}_k}\mu_k=0
\,.
\end{multline*}
In vector form,
\begin{multline*}
  (I-P^T)(m|m)(\alpha_ke^{-\theta^{(m)}_k}\mu_k\,, k\not=m)^T=
\bl(f'_k(e^{\theta^{(m)}_k}-1)\lambda_k\,,k\not=m\br)^T\\+g'_m\bl(
h_m(\theta^{(m)})\br)e^{-\theta^{(m)}_m}\mu_m(p_{mk}\,,k\not=m)^T\,,
\end{multline*}
where $(I-P^T)(m|m)$ represents the matrix that is obtained
 from the matrix $I-P^T$ by deleting
the $m$--th row and the $m$--th column.
Thus, the required condition is satisfied if, entrywise,
\[
  (e^{-\theta^{(m)}_k}\mu_k\,, k\not=m)^T\ge
(I-P^T)(m|m)^{-1}\bl(f'_k(e^{\theta^{(m)}_k}-1)\lambda_k+g'_m(
h_m(\theta^{(m)}))e^{-\theta^{(m)}_m}p_{mk}\mu_m
\,,k\not=m\br)^T\,.
\]
For the special case  of a Jackson network,  the condition reduces to
condition (3.7) in Puhalskii \cite{Puh19}, so,    it is the condition
that
 the facet $F_{\mathcal{K}\setminus\{m\}}$ is
essential. If it is not fulfilled, then there is no optimal trajectory
down the $x_m$--axis, cf. Fig.8 in Puhalskii \cite{Puh19}, where $m=1$\,.

Motion down the $x_m$--axis is  governed by the equation
\begin{multline*}
  \dot q_m(t)=\partial_mH_{\mathcal{K}\setminus \{m\}}(\theta^{(m)})=
f'_m(e^{\theta^{(m)}_m}-1)e^{\theta^{(m)}_m}\lambda_m
+\sum_{l\not=m}
\alpha_l\partial_mh_l(\theta^{(m)})\mu_l\\
+g'_m\bl(
h_m(\theta^{(m)})\br)\partial_mh_m(\theta^{(m)})\mu_m=
f'_m(e^{\theta^{(m)}_m}-1)e^{\theta^{(m)}_m}\lambda_m
+\sum_{l\not=m}
\alpha_le^{\theta^{(m)}_m}e^{-\theta^{(m)}_l}p_{lm}\mu_l\\-g'_m\bl(
h_m(\theta^{(m)})\br)
e^{-\theta^{(m)}_m}
(1-\sum_{l=1}^Kp_{ml}+\sum_{l\not=m}e^{\theta^{(m)}_l}p_{ml})\mu_m
\\=e^{\theta^{(m)}_m}\bl(f'_m(e^{\theta^{(m)}_m}-1)\lambda_m
+\sum_{l\not=m}
\alpha_le^{-\theta^{(m)}_l}p_{lm}\mu_l\br)\\
-g'_m\bl(
h_m(\theta^{(m)})\br)
e^{-\theta^{(m)}_m}
(1-\sum_{l=1}^Kp_{ml}+\sum_{l\not=m}e^{\theta^{(m)}_l}p_{ml})\mu_m
\,.
\end{multline*}
Since $H_{\mathcal{K}\setminus\{m\}}(\theta^{(m)})=0$\,,
 an   argument along the lines of (2.20) in Puhalskii
\cite{Puh19} shows that the momentum $\theta^{(m)}$ is optimal for
going down the $x_m$--axis.
The cost  is $\theta^{(m)}_m$ per unit length.

\def\cprime{$'$} \def\cprime{$'$} \def\cprime{$'$} \def\cprime{$'$}
  \def\cprime{$'$} \def\polhk#1{\setbox0=\hbox{#1}{\ooalign{\hidewidth
  \lower1.5ex\hbox{`}\hidewidth\crcr\unhbox0}}} \def\cprime{$'$}
  \def\cprime{$'$} \def\cprime{$'$} \def\cprime{$'$} \def\cprime{$'$}
  \def\cprime{$'$}

\end{document}